\documentclass[reqno]{amsart}

\usepackage{amsmath, amsfonts, amssymb, hyperref}
\usepackage{graphics}
\usepackage[usenames]{color}

\newcommand{\new}[1]{#1}
\newcommand{\nnew}[1]{#1}
\newcommand{\ok}{ }
\newcommand{\pub}{ }
\newcommand{\no}{ }
{\theoremstyle{plain}
\newtheorem{theorem}{Theorem}[section]
\newtheorem{conjecture}[theorem]{Conjecture}
\newtheorem{proposition}[theorem]{Proposition}
\newtheorem{lemma}[theorem]{Lemma}

}
{\theoremstyle{definition}
\newtheorem{definition}[theorem]{Definition}
\newtheorem*{definition*}{Definition}
\newtheorem{remark}[theorem]{Remark}
\newtheorem{example}[theorem]{Example}
\newtheorem{hypotheses}[theorem]{Hypotheses}
\newtheorem{hypothesis}[theorem]{Hypothesis}
}

\numberwithin{equation}{section}

\newcommand{\gencomp}{\circ}
\newcommand{\nwbody}[1]{\widetilde{#1}}
\newcommand{\syl}[2]{\mathit{Syl}(#1, #2)}
\newcommand{\Wsw}{W_{\mathit{sw}}}
\newcommand{\Wne}{W_{\mathit{ne}}}
\newcommand{\ti}{\times}
\newcommand{\se}[1]{\mathit{se}_{#1}}

\newcommand{\nw}[1]{\mathit{nw}_{#1}}
\newcommand{\nwD}{\nw{D}}
\newcommand{\nbody}[1]{\overline{#1}}
\newcommand{\south}[1]{\mathit{south}_{#1}}
\newcommand{\upbody}[1]{\widehat{#1}}
\newcommand{\rightright}[2]{#1\negthickspace\rightarrow\negthinspace
#2\negmedspace\rightarrow\negmedspace#1\no}
\newcommand{\upup}[2]{#1\uparrow#2\uparrow#1}
\newcommand{\rightup}[2]{#1\negthickspace\rightarrow\negthinspace#2\uparrow#1}
\newcommand{\upright}[2]{#1\uparrow#2\negmedspace\rightarrow\negmedspace#1}

\newcommand{\inv}[1]{\mathrm{inv}(#1)}
\newcommand{\altcomp}{\star}

\begin{document}
\title[Skew Schur classification]{\new{Towards} a combinatorial classification\ok \\ of skew Schur functions}

\author{Peter R. W. McNamara}
\address{Departamento de Matem\'atica, Instituto Superior T\'ecnico, 1049-001 Lisboa, Portugal}
\curraddr{Department of Mathematics, Bucknell University, Lewisburg, PA 17837, USA\no}
\email{\new{peter.mcnamara@bucknell.edu}}

\author{Stephanie van Willigenburg}
\address{Department of Mathematics, University of British Columbia, Vancouver, BC\pub V6T 1Z2\no, Canada}
\email{steph@math.ubc.ca}
\thanks{The second author was supported in part by the National Sciences and Engineering Research Council of Canada.}

\begin{abstract} 
We present a single operation for constructing skew diagrams whose
corresponding skew Schur functions are equal. This combinatorial operation naturally
generalises and unifies all results of this type to date.  Moreover, our operation suggests a
closely related condition that we conjecture is necessary and sufficient
for skew diagrams to yield equal skew Schur functions.
\end{abstract}

\maketitle


\section{Introduction}\label{sec:intro} 
Littlewood-Richardson coefficients arise in a variety of areas of mathematics and therefore not only knowing how to calculate them, but also knowing relations between them, is of importance. More precisely, given partitions $\lambda, \mu, \nu$, the Littlewood-Richardson coefficient $c^\lambda _{\mu\nu}$ arises most prominently in the following three places.  First,\ok in the representation theory of the symmetric group, given  Specht modules $S^\mu$ and $S^\nu$ we have
\begin{equation}\label{equ:LRgeq0}\left(S^\mu \otimes S^\nu\right)\uparrow ^{S_n}=\bigoplus _\lambda c^\lambda _{\mu\nu} S^\lambda .\end{equation}Secondly, considering the cohomology $H^*(Gr(k,n))$ of the Grassmannian, the cup product of Schubert classes $\sigma _\mu$ and $\sigma _\nu$ is given by
$$\sigma _\mu \cup \sigma _\nu = \sum _\lambda c^\lambda _{\mu \nu} \sigma _{\lambda}.$$Lastly, in the algebra of symmetric functions the skew Schur function $s_{\lambda /\mu}$ can be expressed in terms of the basis of Schur functions, $s_\nu$, via
\begin{equation}\label{equ:LRcoeffs}s_{\lambda /\mu}=\sum _\nu c^\lambda _{\mu\nu} s_\nu.\end{equation}

Consequently, knowledge about $c^\lambda _{\mu\nu}$ impacts a number of fields. Examples of knowledge gleaned so far about $c^\lambda _{\mu\nu}$ include a variety of ways to compute them, such as the Littlewood-Richardson rule \cite{LR, Sch77, ThoThesis, Tho78}, inequalities among them that arise from studying eigenvalues of Hermitian matrices \cite{KnutTao}, instances when they evaluate to zero \cite{Purbhoo},  and polynomiality properties that they satisfy 
\cite{DerksenWeyman, KTT04,  Ras04}. 
However, one natural aspect that has yet to be fully exploited 
is that of equivalence classes of equal coefficients. One way to approach this would be to use \eqref{equ:LRcoeffs} and ask when two skew Schur functions are equal. This avenue is worth pursuing since it was recently shown that computing the coefficients $c^\lambda _{\mu\nu}$ is \#P-complete \cite{Narayanan}.

Returning to representation theory, there exist two polynomial representations of $GL_N(\mathbb{C})$ known as Schur modules and Weyl modules. These modules do not form a set of irreducible modules, and so a natural line of enquiry would be to ascertain
when two of them are\no isomorphic. Since these modules are determined up to isomorphism by their characters, we simply need to discover when two characters are\no equal. It so happens that when the modules are indexed by skew diagrams $D$,\ok then the characters are precisely the skew Schur function $s_{D}$ on $N$ variables. In this case we therefore need only determine when two skew Schur functions are\no equal.

The question of skew Schur function equality arises naturally in one other place: the algebra of symmetric functions. As \new{implied} earlier, the skew Schur functions are not a basis for the symmetric functions,\ok and a question currently considered to be intractable is to find all relations among them. In \cite{RSV07} it was shown that the more specific goal of deriving all binomial syzygies between skew Schur functions could be attained by answering the question of equality.  
\new{For this reason and those cited above, we will attempt to classify all skew Schur function equalities.}
In order to do this, we define the following equivalence relation.

\begin{definition*}
For two skew diagrams $D$ and $D'$  we say they are \emph{skew-equivalent} if $s_D=s_{D'}$,\pub and denote this \new{equivalence} by $D\sim D'$.
\end{definition*}

Our question of when two skew Schur functions are equal then reduces to classifying the equivalence classes of $\sim$.

It should be noted that we are not the first to investigate this \new{equivalence}. In \cite{BTV} skew-equivalence was completely characterized for the subset of skew diagrams known as ribbons (or border strips or rim hooks). Their classification involved
\new{a certain composition of ribbons $\alpha$ and $\beta$ that forms a new ribbon
$\alpha \circ 
\beta$.
The idea behind composition operations is that they allow us to construct new equivalences
from equivalences involving smaller skew diagrams.  For example, the results in \cite{BTV} tell
us that if $\alpha \sim \alpha'$ and $\beta \sim \beta'$ are equivalences of ribbons, then $\alpha \circ \beta \sim \alpha' \circ \beta'$.  They were also able to show that the 
size of every equivalence class of $\circ$ is a specific power of $2$.}
The composition $\circ$ was generalised in \cite{RSV07} to include more general skew diagrams $D$ and yielded compositions $\alpha \circ D$ and $D\circ \beta$. A new composition of skew diagrams denoted by $\alpha\circ _\omega D$ for ribbons $\alpha ,\omega$ and skew diagram $D$ was also
introduced, as was the concept of ribbon staircases. These constructions successfully explained almost all skew-equivalences for skew diagrams with up to $18$ cells, but unfortunately $6$ 
skew-equivalences evaded the authors. In this paper we unify all the above constructions into one construction $D\gencomp _W E$ for skew diagrams $D,E$ and $W$. This composition not only provides us with an explanation for all skew-equivalences discovered to date, but
\new{it also suggests necessary and sufficient conditions for skew-equivalence.
Furthermore, it} affords us the possibility to conjecture that 
all equivalence classes are a specific power of $2$ in size. 
More precisely, this paper is structured as follows.

In the next section we review the necessary preliminaries such as skew diagrams and symmetric functions, and recall two identities that will be crucial in our main proof.  The first of these \new{identities} is the 
Hamel-Goulden determinant that was used with much success in \cite{RSV07} to determine 
skew-equivalence, and the second identity is the classical matrix theory result known as Sylvester's Determinantal Identity.  In Section~\ref{sec:gencomp} we describe how to compose two skew diagrams $D$ and $E$ with respect to a third, $W$, to obtain $D\gencomp _W E$.  For ribbons $\alpha, \beta, \omega$ and a skew diagram $D$ we discuss how our composition generalises the composition $\alpha \circ \beta$ of \cite{BTV} and generalises the compositions $\alpha \circ D$, $D\circ \beta$ and $\alpha\circ _\omega D$  plus the notion of ribbon staircases found in \cite{RSV07}. It is also in this section that we state our central theorem, Theorem~\ref{thm:gencomp}. This theorem is the key to proving our sufficient condition for skew-equivalence, and the whole of Section~\ref{sec:proof} is devoted to its proof.  
Finally in Section~\ref{sec:conclusion}, as a consequence of Theorem~\ref{thm:gencomp}, Theorem~\ref{thm:compright} gives our sufficient condition for skew-equivalence.
We propose in Conjecture~\ref{con:iff} that a closely related condition is necessary and sufficient for 
skew-equivalence, and that the size of every equivalence class is a specific power of $2$.
We also derive some conditions under which $D$ is skew-equivalent to its transpose in Proposition~\ref{pro:transpose} and conjecture the that converse is also true.

\subsection*{Acknowledgements} The authors would like to thank
Nantel Bergeron and Vic Reiner for comments that helped to spark productive lines of investigation.
The Littlewood-Richardson calculator \cite{BucSoftware} and the SF package \cite{SteSoftware} aided invaluably in data generation.  \new{We also thank the anonymous referee for thoughtful suggestions that improved the exposition.}

\section{Preliminaries}\label{sec:prelims} 
\subsection{Diagrams}\label{subsec:diagrams}Before we embark on studying skew Schur functions, we need to recall the following combinatorial constructions. We say a \emph{partition}, $\lambda$, of $n$ is a list of positive integers $\lambda _1\geq \lambda _2 \geq \cdots \geq \lambda _k >0$ whose sum is $n$. We denote this by $\lambda \vdash n$, call $k=:\ell(\lambda)$ the \emph{length} of $\lambda$, and call $n$ the \emph{size} of $\lambda$, denoted by $|\lambda |$. For convenience we denote the unique partition of $0$ by $\emptyset$. To every partition $\lambda$ we can associate a subset of $\mathbb{Z}^2$ called a \emph{diagram} that consists of $\lambda _i$ left-justified cells in row $i$. By abuse of notation we also denote this diagram by $\lambda$. In the example below, the symbol $\times$ denotes a cell, although in what follows we may choose to denote cells by numbers, letters or boxes for further clarity.

\begin{example}
$$(3,2,2,1)=\begin{matrix}
\ti&\ti&\ti\\
\ti&\ti\\
\ti&\ti\\
\ti
\end{matrix}.$$
\end{example}
Using this convention for constructing diagrams, we locate cells in the diagram by their row and column indices $(i,j)$, where $i\leq \ell(\lambda)$ and $j\leq \lambda _1$. Moreover, if a cell is contained in row $i$ and column $j$ of a diagram,\ok then we say $c(i,j)=j-i$ is the \emph{content} or \emph{diagonal} of the cell.  We will often use navigational terminology to refer to cells of a diagram.  For example, 
the south (\nnew{respectively} east) border consists of those cells $(i,j)$ such that $(i+1,j)$ (resp.\ $(i,j+1)$)
is not an element of the diagram, 
while the southeast border consists of 
those cells $(i,j)$ such that $(i+1,j+1)$ is not an element of the diagram.
A cell $(i,j)$
is said to be strictly north of a cell $(i',j')$ if $i < i'$, while $(i,j)$ is said to be one position 
northwest of $(i',j')$ if $(i,j)=(i'-1,j'-1)$.

Now consider two diagrams $\lambda$ and $\mu$ such that $\ell (\lambda)\geq \ell(\mu)$ and 
$\lambda _i\geq \mu _i$ for all $i\leq \ell(\mu)$, which we denote by $\mu \subseteq \lambda$. If we locate the cells of $\mu$ in the northwest corner of the set of cells of $\lambda$,\ok then the \emph{skew diagram} $\lambda /\mu$ is the array of cells contained in $\lambda$ but not in $\mu$, where $\lambda /\emptyset = \lambda$. As an example of a skew diagram we have
$$(3,2,2,1)/(2,1)=\begin{matrix}
&&\ti\\
&\ti\\
\ti&\ti\\
\ti
\end{matrix}.$$For convenience we will often refer to generic skew diagrams by capital letters such as $D$. As with partitions we will call the number of cells in $D$ the size of $D$ and denote it by $|D|$. We also consider two skew diagrams to be equal as subsets of the plane if one can be obtained from the other by the addition or deletion of empty rows or columns, or by vertical or horizontal translation.

Any subset of the cells of $D$ that itself forms a skew diagram is said to be a \emph{subdiagram}
of $D$.  If two cells $(i,j)$ and $(i',j')$ satisfy $|i-i'|+|j-j'|=1$,\ok then we say that they are \emph{adjacent},\ok and we similarly say that two \emph{subdiagrams} $D_1$ and $D_2$ are adjacent if \new{$D_1 \cap D_2 = \emptyset$ and} there exists 
a cell in $D_1$ adjacent to a cell in $D_2$.  This concept will play a fundamental role in the pages to follow, but now we will use it to define what it means to be a connected skew diagram. A skew diagram is said to be \emph{connected} if for every cell $d$ with another cell strictly north or east of it, there exists a cell adjacent to $d$ either to the north or to the east. A connected skew diagram is called a \emph{ribbon} (or \emph{border strip} or \emph{rim hook}) if it does not contain the subdiagram $\lambda=(2,2)$.

Given any connected  skew diagram $D$ there exist two natural subdiagrams of $D$, both of which are ribbons. The first is denoted by $\nw{D}$ and is the ribbon that starts at the southwesternmost cell of $D$, traverses the \emph{northwest} border of $D$, and ends at the northeasternmost cell of $D$. The second is denoted by $\se{D}$ and is the ribbon that starts at the southwesternmost cell of $D$, traverses the \emph{southeast} border of $D$, and ends at the northeasternmost cell of $D$. A skew diagram closely related to $\se{D}$ is $\nwbody{D}$, defined in set notation by $\nwbody{D} := D \setminus \se{D}$.

To close this subsection, we recall two symmetries on a skew diagram $D$. The first of these \new{symmetries} is the \emph{transpose} or \emph{conjugate} of $D$, denoted $D^t$, which is obtained by reflecting $D$ along the diagonal that runs from northwest to southeast through all cells with content $0$. The second is the \emph{antipodal rotation} of $D$, denoted $D^*$, which is obtained by rotating $D$ by $180$ degrees in the plane.

\subsection{The algebra of symmetric functions}\label{subsec:Sym}The algebra of symmetric functions has many facets to it, and in this section we review the pertinent details required for our results. More information on this fascinating algebra can be found in \cite{MacD, Sagan, ECII}.

Let $\Lambda ^n$ be the set of all formal power series $\mathbb{Z}[ x_1, x_2, \ldots ]$ in countably many variables that are homogeneous of\new{ }degree $n$ in the $x_i$, and invariant under all permutations of the variables. Then the \emph{algebra of symmetric functions} is
$$\Lambda := \bigoplus _{n\geq 0} \Lambda ^n$$where $\Lambda _0= \mbox{span} \{1\} = \mathbb{Z}$. It transpires that $\Lambda$ is a polynomial algebra in the \emph{complete symmetric functions}, which are defined for all integers $r>0$ by
$$
h_r:= \nnew{\sum_{ 1 \leq i_1 \leq i_2 \leq \cdots \leq i_r }} x_{i_1} x_{i_2} \cdots x_{i_r},$$and $h_0=1$. To obtain a $\mathbb{Z}$-basis for $\Lambda$, let $\lambda=(\lambda_1,\ldots,\lambda_k)$ be a partition of $n$ and let
$$h_\lambda  :=h_{\lambda_1} \cdots h_{\lambda_k},$$for which we find $\{h_\lambda\}_{\lambda \vdash n}$ is a $\mathbb{Z}$-basis for $\Lambda ^n$. However, for the reasons cited in the \no introduction, it is arguable that the most important $\mathbb{Z}$-basis of $\Lambda$ is that consisting of the Schur functions, which we now define as a subset of the skew Schur functions.

Given a skew diagram $D$, we say that $T$ is a \emph{semistandard Young tableau} if $T$ is a filling of the cells of $D$ with positive integers such that:
\begin{itemize}
\item the entries in the rows weakly increase when read from west to east,\ok and
\item the entries in the columns strictly increase when read from north to south.
\end{itemize}
The \emph{skew Schur function} 
$s_D$ is then
\begin{equation}
\label{equ:tableau-Schur-definition}
s_D : = \sum_{T} x^T
\end{equation}
where the sum ranges over all semistandard Young tableaux of shape $D$, and
$$
x^T:= \prod_{(i,j) \in D} x_{T_{ij}}.$$Moreover, the skew Schur function is a \emph{ribbon Schur function} if $D$ is a ribbon,\ok and it\ok is a \emph{Schur function} if for $D=\lambda /\mu$ we have that $\mu =\emptyset$. In this latter case we usually write $s_D=s_\lambda$, which yields another description of $\Lambda$ as  $\Lambda = \bigoplus _{n\geq 0}\ok \Lambda ^n$,\ok where $\Lambda ^n=\mbox{span}\{ s_\lambda | \lambda \vdash n\}$. Using \eqref{equ:tableau-Schur-definition} above, observe that the relationship between the basis of complete symmetric functions and the Schur functions is $h_r=s_{(r)}$, where $(r)$ denotes
the diagram consisting of one row of $r$ cells for $r>0$. It follows from \eqref{equ:LRgeq0} and \eqref{equ:LRcoeffs} that the skew Schur functions are \emph{Schur-positive}, i.e.\ they can be written as a nonnegative linear combination of Schur functions. Moreover, products of Schur functions are also Schur-positive. This follows via the adjointness property with respect to the Hall inner product \cite[Chapter~1, Equation~(5.1)]{MacD}, which says that given partitions $\lambda, \mu, \nu$ with $\mu \subseteq \lambda$
$$
\langle s_\lambda, s_\mu s_\nu \rangle =
\langle s_{\lambda/\mu}, s_\nu \rangle.
$$Hence products of Schur functions and thus products of skew Schur functions are Schur-positive. This property of Schur-positivity will be useful to us later. However, it is a skew Schur function expansion into ribbon Schur functions that we wish to pursue in detail now.

\subsection{Hamel-Goulden determinants}\label{subsec:HGdets}There are a number of ways of expressing a skew Schur function in terms of a matrix determinant involving ribbon Schur functions. For example, the Jacobi-Trudi determinant involves ribbons that are rows or columns, the Giambelli determinant involves hooks, and the Lascoux-Pragacz determinant involves certain more complex ribbons. However, the determinant we will describe generalises all of these, and 
since it was introduced in \cite{HaGo95}, it is known as the Hamel-Goulden determinant.  Our description follows that of \cite{CYY05}.

If $D$ is a skew diagram, then a \new{\emph{decomposition} of $D$ is simply a partition of the elements
of $D$ into disjoint subdiagrams of $D$.}  A \emph{ribbon decomposition} is an ordered decomposition $\Pi=(\theta _1, \ldots, \theta _m)$ of $D$ into ribbons, and furthermore it is an \emph{outside (ribbon) decomposition} if each $\theta _i$ is a ribbon whose southwesternmost (resp.\ northeasternmost) cell lies on either the west or south (resp.\ east or north) border of $D$. In this case we call $\theta _i$ an \emph{outside ribbon}. Given an outside decomposition $\Pi$ of $D$ we say a cell $x\in \theta _i$:
\begin{itemize}
\item \emph{goes north} if the cell adjacent to $x$ to the north is also in $\theta _i$, or $x$ is the northeasternmost cell of $\theta _i$ and lies on the north border of $D$; and
\item \emph{goes east} if the cell adjacent to $x$ to the east is also in $\theta _i$, or $x$ is the northeasternmost cell of $\theta _i$ and lies on the east border of $D$.
\end{itemize}
Observe that every cell in the same diagonal of $D$ either goes north or goes east with respect to $\Pi$. With this in mind the \emph{cutting strip} $\theta (\Pi)$ of $\Pi$ is the unique ribbon occupying the same diagonals as $D$ such that a cell $x\in \theta (\Pi)$ goes north or goes east if and only if the cells in $D$ with content $c(x)$ go north or go east with respect to $\Pi$. Note that each ribbon $\theta _i$ naturally corresponds to the subdiagram of $\theta (\Pi)$ that contains the cells whose contents lie in the interval $[p(\theta _i), q(\theta _i)]$, where $p(\theta _i)$ is the content of the southwesternmost cell of $\theta _i$, and  $q(\theta _i)$ is the content of the northeasternmost cell of $\theta _i$. Extending this notion, we define $\theta[p,q]$ to be the subdiagram of $\theta (\Pi)$ that contains the cells whose contents lie in the interval $[p,q]$ where:
\begin{itemize}
\item $\theta[q+1,q]=\emptyset$, the empty ribbon,\ok and
\item $\theta[p,q]=\mathit{undefined}$ when $p>q+1$.
\end{itemize}
If we define
$$
\theta_i \# \theta_j:=\theta[p(\theta_j),q(\theta_i)],\ok
$$then the \emph{Hamel-Goulden determinant} states that for any outside decomposition $\Pi=(\theta _1, \ldots,  \theta _m)$ of a skew diagram $D$ we have 
\begin{equation}
\label{equ:Hamel-Goulden-formula}
s_D = \det( s_{\theta_i \# \theta_j} )_{i,j=1}^{m}\, ,\ok
\end{equation}
where $s_\emptyset =1$ and $s_{\mathit{undefined}}=0$. We call $( s_{\theta_i \# \theta_j} )_{i,j=1}^{m}$ the \emph{Hamel-Goulden matrix}.

\begin{example} 
If $D=(3,3,3,1)/(1)$,\ok then one possible outside decomposition 
$\Pi=(\theta_1,\theta_2)$ is shown below. The cells in $\theta _i$ are labelled by $i$, and the  cutting strip $\theta(\Pi)$ and  identification of the ribbons $\theta_i$ with intervals of contents within $\theta(\Pi)$ are also shown.
$$
D=
\begin{array}{ccc}
 &1&1\\
1&1&2\\
1&2&2 \\
1&&    
\end{array}
\quad \quad
\theta(\Pi)=
\begin{array}{ccc}
 &\ti&\ti\\
\ti&\ti&\\
\ti&& \\
\ti&& 
\end{array}
\quad\quad
\begin{aligned}
\theta_1&\leftrightarrow \theta[-3,2],\ok \\
\theta_2&\leftrightarrow \theta[-1,1].\ok 
\end{aligned}
$$
The Hamel-Goulden determinant is then
$$
s_D
=\det
\left[
\begin{matrix}
s_{\theta[-3,2]} & s_{\theta[-1,2]} \\
s_{\theta[-3,1]} & s_{\theta[-1,1]} \\
\end{matrix}
\right]
$$
$$
=\det
\left[
\begin{matrix}
  & &  & &\\
 &  
 
 s_{\begin{matrix} &\ti&\ti\\
\ti&\ti&\\
\ti&& \\
\ti&&
    \end{matrix}} 
    
& &

  s_{\begin{matrix}     &\ti&\ti\\
\ti&\ti&\\
&& \\
&&
     \end{matrix}}\\
     
 & & & &\\
 &
 
 s_{\begin{matrix} &\ti&\\
\ti&\ti&\\
\ti&& \\
\ti&&  \end{matrix}} 

& &

  s_{\begin{matrix}    &\ti&\\
\ti&\ti&\\
&& \\
&&

     \end{matrix}}\\
     & & & &\\
\end{matrix}
\right].
$$
\end{example}
Observe that in this case $\theta(\Pi)=\theta _1=\nw{D}$. However, this is not true in general but is a property of the outside decomposition we chose, which we define next.

\begin{definition}\label{def:nw/sedecomp}Given a connected skew diagram $D$ the \emph{southeast decomposition} is an outside decomposition of $D$ that is unique up to reordering. We construct it by choosing the first ribbon to be $\se{D}$. Now consider $D$ with $\se{D}$ removed and iterate the procedure on the remaining skew diagram. If this skew diagram is no longer connected,\ok then iterate on each of the connected components.

We can similarly define the \emph{northwest decomposition} by utilising $\nw{D}$.
\end{definition}

Observe that the above example is a northwest decomposition, and it is straightforward to see that in general with a northwest decomposition, $\Pi$, of a skew diagram $D$, the cutting strip is $\theta (\Pi)=\nw{D}$. Similarly, if we had used a southeast decomposition,\ok then $\theta (\Pi)=\se{D}$. A third outside decomposition that will be useful later is the \emph{horizontal}
or \emph{Jacobi-Trudi decomposition} in which $\theta _i$ is simply row $i$ of the skew diagram $D$. 

The last of our preliminaries is the following result, known as Sylvester's Determinantal Identity, and can be found in
standard matrix theory references, such as \cite{HoJo85}.  This\ok \nnew{identity} will serve a pivotal
role in the proof of Theorem~\ref{thm:gencomp}, where the matrix in question
will be a Hamel-Goulden matrix.  If $M$ is an $n$-by-$n$ matrix, and $A, B \subseteq\{1,\ldots,n\}$, then we let $M[A,B]$ denote the submatrix of $M$ consisting of those entries $(i,j)$ of $M$ having $i \in A$ and $j \in B$.

\begin{theorem}\label{thm:sylvester}
Let $M$ be an $n$-by-$n$ matrix, and let $S \subseteq \{1,\ldots,n\}$.  Define a matrix $\syl{M}{S}$, with
rows and columns indexed by $\{1,\ldots,n\} \setminus S$, by
\[
\syl{M}{S}_{i,j} = \det M[S \cup \{i\}, S \cup \{j\}]
\]
for $i, j \notin S$.  Then
\begin{equation}\label{equ:sylvester}
(\det M) (\det M[S, S])^{n-|S|-1}  = \det \syl{M}{S} . 
\end{equation}
\end{theorem} 

In the case when $S=\{2,\ldots, n-1\}$, \eqref{equ:sylvester} is known as the Desnanot-Jacobi Identity\ok
and is equivalent to Dodgson's condensation formula. 

\section{Compositions of skew diagrams}\label{sec:gencomp}It is now time to recall our equivalence relation that was defined for ribbons in \cite{BTV} and generalised in \cite{RSV07}.

\begin{definition}\label{def:skewequiv} For two skew diagrams $D$ and $D'$  we say they are \emph{skew-equivalent} if $s_D=s_{D'}$, and denote this \new{equivalence} by $D\sim D'$.
\end{definition}

The goal of this paper is to classify skew-equivalence by a condition that is both necessary and sufficient. Fortunately the number of skew-equivalences we need to classify is greatly reduced due to

\begin{proposition} \cite[Section~6]{RSV07}\pub Understanding the equivalence relation $\sim$ on all skew diagrams is equivalent to understanding $\sim$ among \emph{connected} skew diagrams.
\end{proposition}

Consequently, we will \nnew{henceforth assume that} all skew diagrams are connected unless otherwise stated.

We will also make use of the following necessary conditions for skew-equivalence.

\begin{theorem}\label{thm:necconds}
For skew diagrams $D$ and $D'$, if $D \sim D'$\pub then we have:
\begin{enumerate}
\renewcommand{\theenumi}{\roman{enumi}}
\item The number of cells in $D$ equals the number of cells in $D'$.
\item The number of rows in $D$ equals the number of rows in $D'$.
\item $| \nwbody{D}| = |\nwbody{D'}|$.
\end{enumerate} 
\end{theorem}

\begin{proof}
\begin{enumerate}
\renewcommand{\theenumi}{\roman{enumi}}
\item This comes from the definition of $s_D$ in terms of tableaux \eqref{equ:tableau-Schur-definition}.
\item This follows immediately from \cite[Proposition~6.2(ii)]{RSV07},\ok where it was shown that the multisets of row lengths of $D$ and $D'$ are equal.
\item There exists an involution $\omega$ on $\Lambda$ such that $\omega(s_D)=s_{D^t}$ \cite[Chapter~1, Equation~(5.6)]{MacD} and hence $D^t\sim ({D'})^t$. By the previous part this implies that the number of columns in $D$ equals the number of columns in $D'$. Since \new{$| \se{D}|\nnew{+}1$} is the number of rows of $D$ plus the number of columns of $D$, the result now follows by the first part.
\end{enumerate}
\end{proof}

Our approach throughout will be to use known skew-equivalences to construct
skew-equivalences for larger skew diagrams.  
Our basic building blocks will be the skew-equivalences of the following proposition, which is
not hard to prove using the symmetry of $s_D$ and its definition in terms 
of tableaux \eqref{equ:tableau-Schur-definition}.

\begin{proposition}\cite[Exercise~7.56(a)]{ECII}\pub \label{pro:rotation}
For any skew diagram $D$, $D^* \sim D$.
\end{proposition}

The other main ingredient, and the focus of this paper, is a way to put these building blocks
together to construct more complex skew-equivalences.
More specifically, we wish to define
a notion of composition $D \gencomp E$ for skew diagrams $D$ and $E$.  Then
 if $D \sim D'$ and $E \sim E'$, our hope will be that $D \gencomp E \sim D' \gencomp E'$.
Since we wish to generalise and unify the three main operations
of \cite{RSV07}, some care needs to be taken when defining our composition operation,
and some preliminary work is in order.

\begin{definition}
Given skew diagrams $W$ and $E$, we say that $W$ lies in the top (resp.\ bottom)
of $E$ if $W$ appears as a connected subdiagram of 
$E$ that includes the northeasternmost (resp.\ southwesternmost) cell of $E$.

Given two skew diagrams $E_1$ and $E_2$ and a skew diagram $W$ lying
in the top of $E_1$ and the bottom of $E_2$, the \emph{amalgamation of $E_1$ and $E_2$
along $W$}, denoted by $E_1 \amalg_W E_2$, is the new skew diagram obtained from
the disjoint union of $E_1$ and $E_2$ by identifying the copy of $W$ in the top
of $E_1$ with the copy of $W$ in the bottom of $E_2$.
\end{definition}

If $W$ lies in both the top and bottom of $E$, then we will
let $\Wne$ (resp.\ $\Wsw$) denote the copy of $W$ in the top (resp.\ bottom) 
of $E$.
We can also define
\[
E^{\amalg_W m} = \underbrace{E \amalg_W E \amalg_W \cdots \amalg_W E}
_{m \mbox{\scriptsize{ factors}}}
:= \nnew{(\cdots(( E \amalg_W E)\amalg_W E) \amalg_W \cdots ) \amalg_W E}.
\]

\begin{example}\label{exa:amalg} The skew diagram $E$ given by
\[
E = 
\begin{matrix}
& & & & & & \ti & \ti \\
& & & \ti & \ti & \ti & \ti & \ti \\
& \ti & \ti & \ti & \ti \\
\ti & \ti & \ti & \ti 
\end{matrix}
\]
has
\[
W = \begin{matrix} & \ti & \ti \\ \ti & \ti & \ti \end{matrix}
\]
lying in its top and bottom.  We see that
\[
E \amalg_W E = 
\begin{matrix}
& & & & & & & & & & & w & w \\
& & & & & & & & \ti & \ti & w & w & w \\
& & & & & & w & w & \ti & \ti \\
& & & \ti & \ti & w & w & w & \ti \\
& w & w & \ti & \ti \\
w & w & w & \ti
\end{matrix},
\]
where we use the symbol $w$ to denote the cells of copies of $W$.
Notice that $V = \begin{smallmatrix} & \ti & \ti \\ \ti & \ti \end{smallmatrix}$ also
lies in the top and bottom of $E$, and that $E \amalg_{V} E$ is the same skew diagram
as $E \amalg_W E$.
\end{example}

\begin{example}\label{exa:wempty}
For complete generality, we will also say that when $W = \emptyset$, $W$ lies in the top and
bottom of any skew diagram $E$.  In this case, we\ok will identify $\Wsw$ with the west edge
of the southwesternmost cell of $E$.  Similarly, we will identify $\Wne$ with the 
east edge of the northeasternmost cell of $E$.   
For example, if 
$E = (3,3,2)/(1)$, then $\Wsw$ and $\Wne$ would be identified with the thicker edges as shown:\ok
\setlength{\unitlength}{4mm}
\[
\begin{picture}(3,3)(0,0)
\put(0,0){\line(1,0){2}}
\put(0,1){\line(1,0){3}}
\put(0,2){\line(1,0){3}}
\put(1,3){\line(1,0){2}}
\multiput(1,0)(1,0){2}{\line(0,1){3}}
\multiput(0,1)(3,0){2}{\line(0,1){1}}
\thicklines
\multiput(-0.02,0)(0.02,0){3}{\line(0,1){1}}
\multiput(2.98,2)(0.02,0){3}{\line(0,1){1}}
\put(3.5,1.5){\nnew{.}}
\end{picture}
\]
Then $E \amalg_\emptyset E$ is the skew diagram $(6,6,5,3,2)/(4,3,1)$.  
\end{example}

Now is a good time to introduce some assumptions on $E$ and $W$ that we will
need for our results to hold.

\begin{hypotheses}\label{hyp:initial}
Suppose that $E$ is a skew diagram having $W$ lying in its top and bottom.
We assume that $E$ and $W$ satisfy the following conditions:
\begin{enumerate}
\renewcommand{\theenumi}{\Roman{enumi}}
\item \label{ite:maximal}
$W$ is maximal in the following sense: 
there does not exist a skew diagram $W' \supsetneq W$ that occupies the same set of diagonals
as $W$ and that also lies in the top and bottom of $E$.
\item \label{ite:separation}
$\Wne$ and $\Wsw$ are separated by at least one diagonal.
In other words,
there is at least one diagonal
between $\Wne$ and $\Wsw$ that intersects neither $\Wne$ nor $\Wsw$. 
\item \label{ite:complement}
The complement in $E$ of either copy of $W$ is a connected skew diagram.
\end{enumerate}
\end{hypotheses}

\begin{remark}
Analogues of Hypotheses~\ref{ite:separation} and \ref{ite:complement}
are also necessary for the results in \cite[Section~7.2]{RSV07}.
Notice that the $V$ of Example~\ref{exa:amalg} fails to satisfy Hypothesis~\ref{ite:maximal}.
\new{(It also fails to satisfy Hypothesis~\ref{ite:complement}.)
As we saw, however, $E \amalg_{V} E = E \amalg_W E$, and it will be apparent from our
definition of the composition operation that
$D \gencomp E$ for any $D$ is the same whether we work
with $V$ or $W$. Therefore, we lose no generality
when we impose Hypothesis~\ref{ite:maximal}---it will just make the statements
of some of our results simpler.  }
\end{remark}

Hypotheses~\ref{ite:separation} and \ref{ite:complement}
tell us much about the structure of $E$.  Let $O$ denote the subdiagram of $E$ that
results when we delete both copies of $W$.  
We will write
$E=WOW$ to mean that 
$W$ lies in the top and bottom of $E$ and that $O$ is the subdiagram of $E$
that results when we delete both copies of $W$.  Since $E$ is assumed to be connected,
Hypotheses~\ref{ite:separation} and \ref{ite:complement} tell us that $O$ is a non-empty\no
connected skew diagram.  

Let us say that the lower (resp.\ upper) copy of $W$ is \emph{horizontally attached
to $O$} if the southwesternmost (resp.\ northeasternmost) cell of $O$ has a cell of $W$ one 
position to its west (resp.\ east).  Similarly, we say that the lower (resp.\ upper) copy of $W$ is 
\emph{vertically attached to $O$} if the southwesternmost (resp.\ northeasternmost) cell of $O$ has
a cell of $W$ one position to its south (resp.\ north).  Since 
$W$ is a skew diagram and $E$ is connected, 
each copy of $W$ in $E$ is either horizontally or vertically
attached to $O$, but not both.  
Therefore, we are in one
of the following four cases:
\begin{enumerate}
\renewcommand{\theenumi}{\alph{enumi}}
\item Both copies of $W$ are horizontally attached to $O$, written $E=\rightright{W}{O}$.
\item Both copies of $W$ are vertically attached to $O$, written $E=\upup{W}{O}$.
\item The lower copy of $W$ is horizontally attached to $O$, while the upper
copy of $W$ is vertically attached to $O$, written $E=\rightup{W}{O}$.
\item The lower copy of $W$ is vertically attached to $O$, while the upper
copy of $W$ is horizontally attached to $O$, written $E=\upright{W}{O}$.
\end{enumerate}

We are almost ready to define\no composition of general skew diagrams.  One issue that lengthens
the definition of the composition of $D$ and $E$ with respect to $W$ is that the definition 
varies according to the cases (a), (b), (c) and (d) above.  As justification for this variation, 
consider the following
diagrams that can be created, starting with two copies $E_1$ and $E_2$ of $E$:
\begin{enumerate}
\renewcommand{\theenumi}{\Alph{enumi}}
\item Position $E_2$ so that the lower copy of $W$ in $E_2$ is one position northwest
of the upper copy of $W$ in $E_1$.
\item Position $E_2$ so that the lower copy of $W$ in $E_2$ is one position southeast
of the upper copy of $W$ in $E_1$.
\item Form $E_1 \amalg_W E_2$ and translate an extra copy of $W$ one position
southeast from $E_1 \cap E_2$.
\item Form $E_1 \amalg_W E_2$ and translate an extra copy of $W$ one position
northwest from $E_1 \cap E_2$.
\end{enumerate}

\begin{figure}[htpb]
\[
\scalebox{.6}{\includegraphics{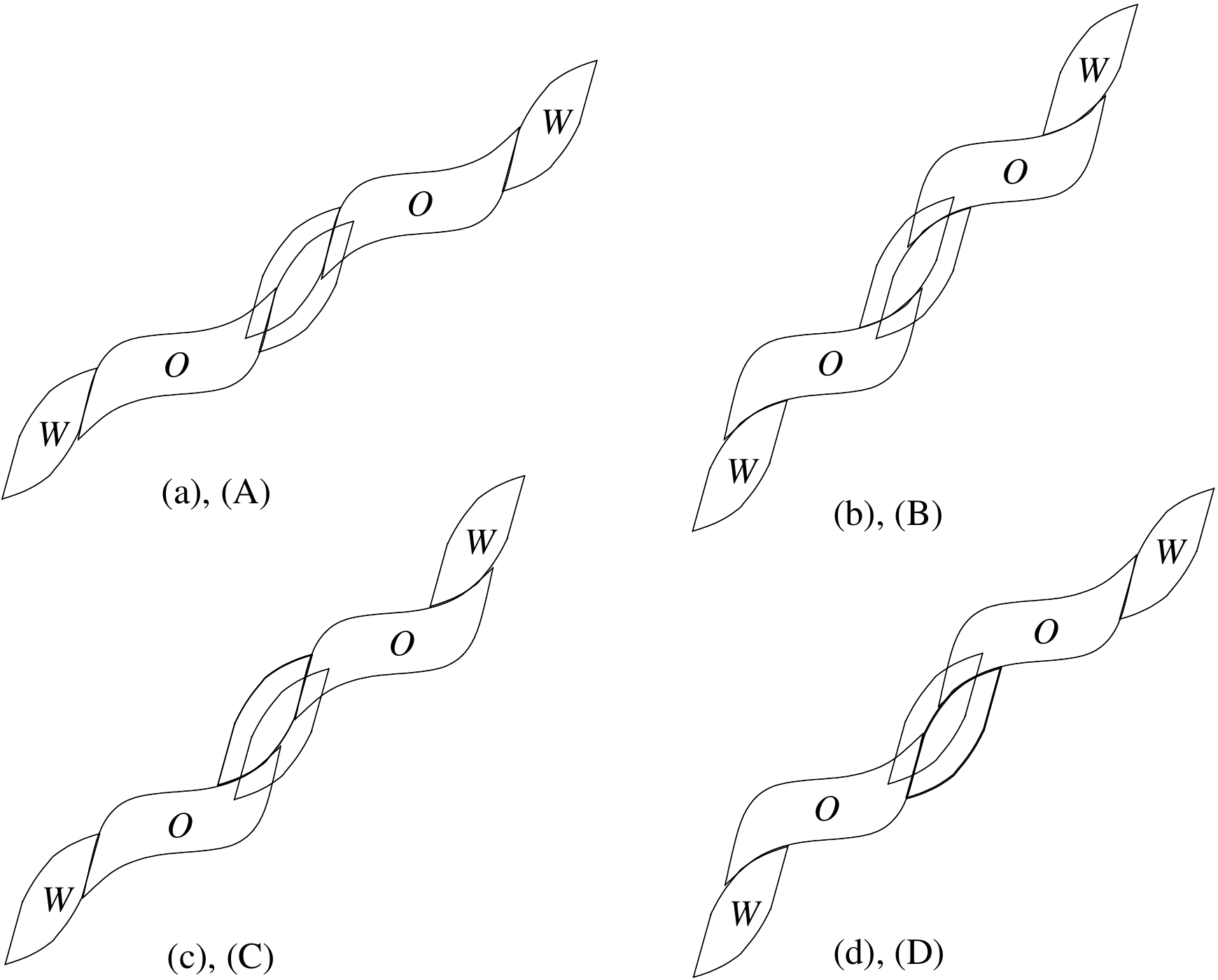}}
\]
\caption{$E_1 \cdot_W E_2$ in the four cases}
\label{fig:cdotoperation}
\end{figure}

The key observation is that in each of the four cases (a), (b), (c) and (d), 
exactly one of these four diagrams
is a skew diagram, namely the diagram with the corresponding letter label.  
This observation effectively consists of sixteen assertions, and we leave
their checking as an exercise for the reader that will reinforce the ideas introduced so far.  
In each of the four cases (a), (b), (c) and (d), we let $E_1 \cdot_W E_2$ denote the 
skew diagram constructed in (A), (B), (C) and (D) respectively.  See Figure~\ref{fig:cdotoperation}
for an illustration.  

\begin{remark}
We see that there is a fundamental difference between the set-up for the cases
$\rightright{W}{O}$, $\upup{W}{O}$ and the cases $\rightup{W}{O}$, $\upright{W}{O}$.
In a certain sense, this \new{difference} is to be expected, since it turns out that $\rightright{W}{O}$
and $\upup{W}{O}$ are involved in generalising
the composition and amalgamated composition operations of \cite{RSV07}, while
$\rightup{W}{O}$ and $\upright{W}{O}$ are involved in generalising the ribbon staircase operation.
The real strength of our framework will be highlighted by the statements of the results that
follow, where all four cases can be treated as one.  
\end{remark}

We are finally ready to define the composition of general skew diagrams.

\begin{definition}\label{def:gencomp}
For skew diagrams $D$, $E$ with $E=WOW$ subject to Hypotheses~\ref{hyp:initial}, 
we define the composition $D \gencomp_W E$ with respect to $W$ as follows.  Every cell
$d$ of $D$ will contribute a copy of $E$, denoted $E_d$, in the plane. 
The set of copies $\left\{E_d\ |\ d \in D\right\}$ are combined according to the following
rules:
\begin{enumerate}
\item[(a), (b)] Suppose $E = \rightright{W}{O}$ or $E=\upup{W}{O}$.
\begin{enumerate}
\renewcommand{\theenumii}{\roman{enumii}}
\item If $d$ is one position west of $d'$  in $D$, then $E_d$ and $E_{d'}$ appear
in the form $E_d \amalg_W E_{d'}$.
\item If $d$ is one position south of $d'$ in $D$, then $E_d$ and $E_{d'}$ appear
in the form $E_d \cdot_W E_{d'}$.
\end{enumerate}
\item[(c), (d)] If $E=\rightup{W}{O}$,\ok then we consider the northwest ribbon decomposition
of $D$, while if $E=\upright{W}{O}$,\ok then we consider the southeast ribbon decomposition of 
$D$.  
\begin{enumerate}
\renewcommand{\theenumii}{\roman{enumii}}
\item If $d$ is one position west of $d'$ on the same ribbon in $D$, then $E_d$ and $E_{d'}$
appear in the form $E_d \amalg_W E_{d'}$.  
\item If $d$ is one position south of $d'$ on the same ribbon in $D$, then $E_d$ and $E_{d'}$
appear in the form $E_d \cdot_W E_{d'}$.  
\item If $d$ is one position southeast of $d'$ in $D$, then $E_d$ appears one position 
southeast of $E_{d'}$. 
\end{enumerate}
\end{enumerate}
Additionally, we will use the convention that $\emptyset \gencomp_W E = W$ and 
that $D \gencomp_W E$ is undefined when $D$ is undefined.
\end{definition}

\begin{example}\label{exa:gencompdef}
Identifying the cells of $D$ with integers, and labelling the cells of the copies of
$W$ in $E$ with the letter $w$, suppose
\[
D = \begin{matrix} 1 & 2 \\ 3 & 4 \end{matrix} \mbox{\ \ and\ \ } 
E = \begin{matrix}  &&&w \\ w&\ti&\ti&w \\ w&\ti \end{matrix}\ .
\]
Then $E=\rightright{W}{O}$,\ok and so $D \gencomp_W E$ is the skew diagram
\[
\begin{matrix}
&&&&&&&& 2 \\
&&&&& \ti & 2 & 2 & 2 \\
&& 1 & 1 & 1 & \ti & \ti \\
&& 1 & \ti & 4 & 4 & 4 \\
3 & 3 & 3 & \ti & 4 \\
3 & 3 \\
\end{matrix},\ok
\]
where a cell is labelled by $\ti$ if it is an element of $E_d$ for more than one $d \in D$,
and otherwise is labelled by $d$ when it is an element of $E_d$.  

Alternatively, if 
\[
E = \begin{matrix}  &&w \\ &&w \\ w&\ti&\ti \\ w&\ti \end{matrix}\ ,
\]
then $E = \rightup{W}{O}$,\ok and so $D \gencomp_W E$ is the skew diagram
\[
\begin{matrix}
&&&&&& 2 \\
&&&&&& 2 \\
&&&& \ti & 2 & 2 \\
&&&& \ti & \ti \\
&& \ti & 1 & 1 & 4 \\
&& \ti & \otimes & 4 & 4 \\
3 & 3 & 3 & \otimes & 4 \\
3 & 3 
\end{matrix},
\]
where $\otimes$ denotes an element of both $E_4$ and $E_3 \cdot_W E_1$.  
\end{example}

\begin{example}\label{exa:obtainingcirc}
If $W$ is empty, then referring to Example~\ref{exa:wempty}, it is natural to 
consider $E$ to be of the form $\rightright{W}{O}$.  If at least one of $D$ and $E$
is a ribbon, $\gencomp_\emptyset$ becomes the composition denoted simply by $\circ$
in \cite{RSV07}.  When both $D$ and $E$ are ribbons, $D \gencomp_\emptyset E$
also corresponds to $D \circ E$ of \cite{BTV}.  When neither $D$ nor $E$ is a ribbon,
{$\gencomp_\emptyset$} behaves like $\circ$, except that we allow overlaps to occur
among copies of $E$.  
To see this
in action, take $D$ as in the previous example,
and let  $E=\begin{smallmatrix} \ti & \ti \\ \ti & \ti \end{smallmatrix}$.
Then 
\[
D \gencomp_{\emptyset} E = 
\begin{matrix}
&&& 2 & 2 \\
& 1 & 1 & 2 & 2 \\
& 1 & \ti & 4 \\
3 & 3 & 4 & 4 \\
3 & 3
\end{matrix}.
\]
We note in passing that this is the same skew diagram that appears in \cite[Remark~7.10]{RSV07},
and was the first motivating example for the work of the current article.
\end{example}

\begin{example}
The previous example demonstrates how one of the three main operations of \cite{RSV07}
is obtained as a special case of our composition operation.
The other two operations are also obtained as special cases:
the amalgamated composition operation corresponds to certain cases of 
$D \gencomp_W E$ with $D$ and $W$ non-empty\no ribbons and 
$E$ of the \nnew{form} $\rightright{W}{O}$ or $\upup{W}{O}$.
On the other hand, if $E$ is a ribbon of the form $\rightup{W}{O}$ or $\upright{W}{O}$, 
then $D \gencomp_W E$ is a ribbon staircase construction.
\end{example}

\begin{remark}\label{rem:horizdecomp}
We see in Definition~\ref{def:gencomp} that $E=\rightup{W}{O}$ is associated with the northwest ribbon decomposition of $D$,
while $E=\upright{W}{O}$ is associated with the southeast ribbon decomposition.
We remark that $E=\rightright{W}{O}$ and $E=\upup{W}{O}$ should both be
associated with the horizontal ribbon decomposition.  Indeed, (i) above for the (a), (b) case
could equivalently state that if
$d$ is one position west of $d'$ on the \new{\emph{same ribbon}} of the horizontal decomposition
of $D$, then $E_d$ and $E_{d'}$ appear in the form $E_d \amalg_W E_{d'}$.  This is
the same rule as (i) for the (c), (d) case.
More importantly, our proof of Theorem~\ref{thm:gencomp} for 
$E=\rightright{W}{O}$ and $E=\upup{W}{O}$ will exploit the horizontal decomposition
of $D$, in the same way that our proof for $E=\rightup{W}{O}$ and $E=\upright{W}{O}$
will exploit the northwest and southeast
decompositions.  
\end{remark}

When $E=\rightup{W}{O}$, we have defined $D \gencomp_W E$ in terms of the
northwest decomposition of $D$.  Before proceeding, we give an alternative
definition of $D \gencomp_W E$, this time in terms of the southeast decomposition 
of $D$.  This is a necessary tool for proving Lemma~\ref{lem:gencomp}(\ref{ite:transpose})
below.  The reader that is focussing solely on the statements of the main results can safely skip
the next definition, example and lemma. 

First, we will give an alternative definition of $D \gencomp_W E$ that
highlights the way in which its structure is affected by the extra copies of $W$
that arise from copies of $E$ that are in the form $E_d \cdot_W E_{d'}$.  
Let $W_1$ denote the extra copy of $W$ that is added to $E_d \amalg_W E_{d'}$
to form $E_d \cdot_W E_{d'}$.  We observe that $W_1$
will be covered by a copy $E_{d''}$ for some $d'' \in D$ in many
situations.  Specifically, if $d'$ has a cell $d''$ one position to its southeast, then
$W_1$ will be covered by $E_{d''}$.  Therefore, we could also define 
$D \gencomp_W E$ as follows.

\begin{definition}\label{def:gencomp2}
For skew diagrams $D$ and $E$ with $E$ of the form $\rightup{W}{O}$, we define
$D \gencomp_W E$ as follows.  Every cell $d$ of $D$ will contribute a copy of $E$, 
denoted $E_d$, in the plane.  
Considering the northwest ribbon decomposition of $D$, we position the 
set of copies $\{E_d\ |\ d \in D\}$ according to the following rules:
\begin{enumerate}
\renewcommand{\theenumi}{\roman{enumi}}
\item If $d$ is one position west or south of $d'$ on the same ribbon in $D$, then $E_d$ and $E_{d'}$
appear in the form $E_d \amalg_W E_{d'}$.  
\item If $d$ is one position southeast of $d'$ in $D$, then $E_d$ appears one position 
southeast of $E_{d'}$. 
\end{enumerate}
Furthermore, if $d$ is one position south of $d'$ on the same ribbon in $D$ and both $d$ and $d'$ are
elements of $\se{D}$, then add a copy of $W$ one position southeast of  $E_d \cap E_{d'}$. 
\end{definition}

\begin{example}\label{exa:extraws}
Suppose 
\[
D = 
\begin{matrix}
& a_5 \\
a_3 & a_4 \\
a_2 & b \\
a_1 
\end{matrix}
\mbox{\ \ and\ \ }
E =
\begin{matrix}
&&& w & w \\
&&& w & w \\
w & w & \ti & \ti \\
w & w & \ti
\end{matrix}.
\]
Then a representative diagram of 
$D \gencomp_W E$ is shown in Figure~\ref{fig:extraws}.  The normal lines represent
the contribution from the $a_i$ ribbon,\pub the dotted lines represent $E_b$,\ok while the 
bold lines represent the extra copies of $W$.
\begin{figure}[tpb]
\[
\scalebox{.6}{\includegraphics{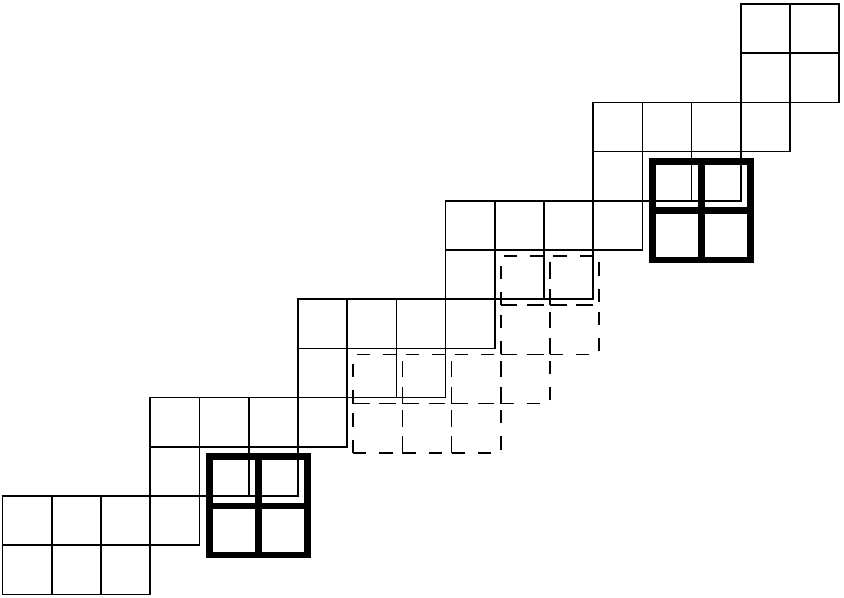}}
\]
\caption{The real contribution from extra copies of $W$}
\label{fig:extraws}
\end{figure}
\end{example}

When $E=\rightup{W}{O}$, we are now in a position to give a definition of $D \gencomp_W E$ 
in terms of the southeast decomposition of $D$, as promised.  It should be compared
with the relevant part of Definition~\ref{def:gencomp}.

\begin{lemma}\label{lem:southeastdef}
For skew diagrams $D$ and $E$ with $E$ of the form $\rightup{W}{O}$, suppose we define a
diagram $D \altcomp_W E$ as follows.  Every cell
$d$ of $D$ will contribute a copy of $E$, denoted $E_d$, in the plane. 
Considering the southeast decomposition of $D$, we position the set of copies 
$\left\{E_d\ |\ d \in D\right\}$ according to the following rules:
\begin{enumerate}
\renewcommand{\theenumi}{\roman{enumi}}
\item If $d$ is one position west of $d'$ on the same ribbon in $D$, then $E_d$ and $E_{d'}$
appear in the form $E_d \amalg_W E_{d'}$.  
\item If $d$ is one position south of $d'$ on the same ribbon in $D$, then $E_d$ and $E_{d'}$
appear in the form $E_d \cdot_W E_{d'}$.  
\item If $d$ is one position northwest of $d'$ in $D$, then $E_d$ appears one position 
southeast of $E_{d'}$. 
\end{enumerate}
Then $D \altcomp_W E = D \gencomp_W E$.
\end{lemma}

\begin{proof}
In the same way that we extracted Definition~\ref{def:gencomp2} from 
Definition~\ref{def:gencomp}, we could extract an analogue of Definition~\ref{def:gencomp2}
from the definition of $D \altcomp_W E$.  This analogue states that,
considering the southeast ribbon decomposition of $D$, we construct 
$D \altcomp_W E$ by positioning the 
set of copies $\{E_d\ |\ d \in D\}$ according to the following rules:
\begin{enumerate}
\renewcommand{\theenumi}{\roman{enumi}}
\item If $d$ is one position west or south of $d'$ on the same ribbon in $D$, then $E_d$ and $E_{d'}$
appear in the form $E_d \amalg_W E_{d'}$.  
\item If $d$ is one position northwest of $d'$ in $D$, then $E_d$ appears one position 
southeast of $E_{d'}$. 
\end{enumerate}
Furthermore, if $d$ is one position south of $d'$ on the same ribbon in $D$ and both $d$ and $d'$ are
elements of $\nw{D}$, then add a copy of $W$ one position southeast of $E_d \cap E_{d'}$. 

Thinking in terms of Definition~\ref{def:gencomp2} and this analogue, 
it is easy to see that an outside ribbon 
of  size $l$ in $D$ contributes $E^{\amalg_W l}$ to $D \gencomp_W E$ \new{and $D \altcomp_W E$.}
Our proof now divides into two parts.  We will first show that these contributions
from the ribbons give the same result, whether we work with the northwest decomposition, 
as in $D \gencomp_W E$, or with the southeast decomposition, as in $D \altcomp_W E$.
To finish, we show that the effect of the extra copies of $W$ from the final sentence of
Definition~\ref{def:gencomp2} on $D \gencomp_W E$ is the same as the effect of the extra 
copies of $W$ from the analogue on $D \altcomp_W E$.

It is clear that $\nwD$ and $\se{D}$ have the same size $l_1$, and \new{so they will contribute 
the same $E_1 = E^{\amalg_W l_1}$ to $D \gencomp_W E$ and $D \altcomp_W E$ respectively.
Furthermore, we can see that $D \setminus \nwD$ and $D \setminus \se{D}$ are the same skew diagram $D_1$.  It will not affect our argument that
$D_1$ need not be connected.  Repeating this process, we could
remove a maximal size ribbon $r_{\mathit{nw}}$,
that is an element of the northwest decomposition of $D$, 
from the northwest border of some connected component $D'_1$ of $D_1$.  We could also remove a maximal size
ribbon $r_{\mathit{se}}$, that is an element of the southeast decomposition of $D$, 
from $D'_1$.
Again, removing either $r_{\mathit{nw}}$ and $r_{\mathit{se}}$ from $D'_1$ results in the same skew diagram $D_2$.
The key point is that
$r_{\mathit{nw}}$ and $r_{\mathit{se}}$ have the same size $l_2$ so, 
while they aren't in general the same ribbon, they contribute the same
$E_2=E^{\amalg_W l_2}$ to $D \gencomp_W E$ and $D \altcomp_W E$ respectively.  Furthermore, since 
$r_{\mathit{nw}}$ and $r_{\mathit{se}}$ have the same set of contents, the copies of $E_2$ will
have the same positions relative to the copies of $E_1$.}

Repeating this process, we will see that the contributions of the form $E^{\amalg_W l}$
to $D \gencomp_W E$ from the ribbons of the northwest decomposition
will be the same,\no and have the
same relative positions, as the analogous contributions to $D \altcomp_W E$ from the ribbons
of the southeast decomposition.  

It remains to show that the extra copies of $W$ will be the same in both cases.
This amounts to showing that the following are equivalent:
\begin{enumerate}
\item There exists $d$ with content $i$ one position south of
$d'$ on the same ribbon of the northwest decomposition of $D$, and both $d$ and
$d'$ are elements of $\se{D}$.
\item There exists $c$ with content $i$ one position south of
$c'$ on the same ribbon of the southeast decomposition of $D$, and both $c$ and
$c'$ are elements of $\nwD$.
\end{enumerate}
Suppose (1).  Since $D$ is a skew diagram, and since $d$ and $d'$ are on the same ribbon
of the northwest decomposition of $D$, we see that $d \in \nwD$ if and only if $d' \in \nwD$.
If $d, d' \in \nwD$, then we conclude (2), since $d$ and $d'$ are on the same ribbon of the
southeast decomposition of $D$.  If $d, d' \not\in \nwD$, then there exist cells
$e$ and $e'$ one position northwest of $d$ and $d'$ respectively.  We see that since $d$ and $d'$
are on the same ribbon in both the northwest and southeast decompositions of $D$, the 
same must apply to $e$ and $e'$.  For the same reasons as before, $e \in \nwD$ if and only if
$e' \in \nwD$.  Repeating this process, working in a northwesterly direction, we will eventually
arrive at $c$ and $c'$ with the required properties.  Since the process is clearly reversible,
the result follows.  
\end{proof}

Readers may wish to check their understanding of the definition of $D \gencomp_W E$
by filling in the details in the proof of the following lemma.  Part (i) is of obvious importance,
while (ii) and (iii) will save us much effort in the proof of Theorem~\ref{thm:gencomp}.

\begin{lemma}\label{lem:gencomp}
For skew diagrams $D$ and $E = WOW$ we have: 
\begin{enumerate}
\renewcommand{\theenumi}{\roman{enumi}}
\item $D \gencomp_W E$ is a skew diagram.
\item \label{ite:star} $ (D \gencomp_W E)^* = D^* \gencomp_{W^*} E^* .$
\item \label{ite:transpose} \nnew{If $W \neq \emptyset$, then} $ (D \gencomp_W E)^t = D^* \gencomp_{W^t} E^t .$
\end{enumerate}
\end{lemma}

\begin{proof}
\begin{enumerate}
\renewcommand{\theenumi}{\roman{enumi}}
\item We can prove this by induction on the number $r$ of ribbons in the appropriate ribbon
decomposition of $D$.  
Suppose $E=\rightright{W}{O}$ or $E=\upup{W}{O}$.  
If $r=1$, then $D \gencomp_W E$ is just $E^{\amalg_W |D|}$.  For $r > 1$, the result follows
from the fact that both $E \cdot_W E$ and $D$ are skew diagrams.  
Now suppose $E=\rightup{W}{O}$ or $E=\upright{W}{O}$.  If $r=1$, then
$D \gencomp_W E$ is $E^{\amalg_W |D|}$ with some copies of $W$ added as appropriate.
The result will be a skew diagram since $E \cdot_W E$ is a skew diagram.  Now suppose 
$r>1$.  First observe how the ribbons of the northwest/southeast 
decomposition of $D$ ``nest'' in each other.  We can check that 
the contributions to $D \gencomp_W E$ of the various
ribbons of the northwest/southeast decomposition of $D$ nest with each other in an 
analogous way.  In the same way that $E \cdot_W E$ is a skew diagram, we conclude
that $D \gencomp_W E$ is a skew diagram.

\item We omit the proof since this identity is straightforward to check using 
Definition~\ref{def:gencomp}.
Observe, though, that if $E = \rightup{W}{O}$, then 
$E^* = \upright{W^*}{O^*}$.

\item Notice that if $E = \rightright{W}{O}$, then $E^t = \upup{W^t}{O^t}$.  We omit
the proof of the cases $E = \rightright{W}{O}$ and $E=\upup{W}{O}$, since
then the identity is checked easily, as in (\ref{ite:star}).  We next consider the
case when $E=\rightup{W}{O}$. 
The ribbons of the northwest decomposition of 
$D$ correspond exactly to the ribbons of the southeast
decomposition of $D^*$.  Comparing Definition~\ref{def:gencomp} in the case $E=\rightup{W}{O}$
with the definition of $D \gencomp_W E$ from the statement of 
Lemma~\ref{lem:southeastdef} yields the result.  
If $E=\upright{W}{O}$,\ok then the proof is similar, since we could develop an appropriate analogue of 
Lemma~\ref{lem:southeastdef}. 
\end{enumerate}\end{proof}

We can now work directly towards the statement of Theorem~\ref{thm:gencomp}, which
expresses $s_{D \gencomp_W E}$ in terms of $s_D$ and $s_E$, and thus
serves as the foundation for all our skew-equivalence proofs.  As mentioned in 
Example~\ref{exa:obtainingcirc}, a \new{feature of our definition of $D \gencomp_W E$
is that we allow overlaps among the copies of $E$, whereas composition operations in prior work 
do not.}  At some point, we must
obtain an understanding of, and account for, these overlaps.
This motivates the following definition of
the skew diagrams $\nbody{W}$ and $\nbody{O}$.

\begin{definition}\label{def:nbody}
Consider the infinite skew diagram
\begin{equation}\label{equ:infiniteamalg}
\mbox{\boldmath{$E$}} := E^{\amalg_W \infty} = \cdots \amalg_{W} E \amalg_{W} E \amalg_{W} \cdots . 
\end{equation}
For every copy $O_1$ of $O$ in $\mbox{\boldmath{$E$}}$ we define
\[
\nbody{O_1} = \{ (i,j) \in O_1\ |\ (i+1,j+1)\in O_1\}.
\]
For every copy $W_1$ of $W$ we define
\[
\nbody{W_1} =  \{ (i,j) \in \mbox{\boldmath{$E$}}\ |\ (i+1,j+1) \in W_1\}   \cup 
\{(i,j) \in W_1\ |\ (i+1, j+1) \in \mbox{\boldmath{$E$}}\}.
\]
Clearly, every copy $O_1$ of $O$ defines the same diagram $\nbody{O_1}$, which
we denote simply by $\nbody{O}$. Similarly, we define $\nbody{W}$.
\end{definition}

\begin{example}\label{exa:nbody}
If 
\[
E = 
\begin{array}{cccccccc}
& & & & & & \bar{w} & \bar{w} \\
& & & \bar{\ti} & \ti & w & w & w \\
&  \bar{w} & \bar{w} & \ti & \ti \\
w & w & w & \ti
\end{array},
\]
then the four cells labelled $\bar{w}$ denote two copies of $\nbody{W}$, 
while the single $\bar{\ti}$
denotes $\nbody{O}$.  In general, however, $\nbody{W}$ need not be a subset of $W$.
For example, we can have
\[
E = \begin{array}{ccc}
& \bar{\ti} & \bar{w} \\
\bar{w} & \ti & w \\
w & \ti
\end{array},
\]
where the cells on the top row comprise one copy of $\nbody{W}$.  
Part of a second copy of $\nbody{W}$\no is also shown. 
\end{example}

Let us make some observations about Definition~\ref{def:nbody}:  
\begin{itemize}
\item $\nbody{O_1}$ is nothing more than $\nwbody{O_1}$.  However, 
using the notation $\nwbody{O}$ in what follows would have the potential to 
cause confusion, and using $\nbody{O}$ instead will help to keep
our notation consistent.
\item We see that
$\nbody{O}$ and $\nbody{W}$ are the shapes that result when we remove the infinite 
southeast ribbon from {\boldmath{$E$}}.  
\item \new{ }Observe that $\nbody{O}$ and $\nbody{W}$ are skew diagrams and that neither
one need be connected.
\end{itemize}

It turns out that we will need one further assumption about the structure of $E$.
We conjecture below that this final assumption encompasses exactly what we
need for our expression for $s_{D \gencomp_W E}$ to hold. 
In $\mbox{\boldmath{$E$}}= E^{\amalg_W \infty}$,  Hypothesis~\ref{ite:separation} tells us
that no two copies of $\nbody{W}$ will be adjacent.  

\begin{hypothesis}\label{hyp:crucial}
Suppose that $E=WOW$.  Assume that $E$ satisfies the following condition:
\begin{enumerate}
\addtocounter{enumi}{3} \renewcommand{\theenumi}{\Roman{enumi}}
\item \label{ite:crucial}
In $\mbox{\boldmath{$E$}}$, no copy of $\nbody{O}$ is adjacent to a copy of $\nbody{W}$.
\end{enumerate}
\end{hypothesis}

\begin{remark}
\new{
Suppose we construct a second
copy of {\boldmath{$E$}} which is the translation of {\boldmath{$E$}} one position
northwest.  Then $\nbody{O}$ and $\nbody{W}$ are exactly the shapes that form
the overlap of the two copies of {\boldmath{$E$}}.
Intuitively,  Hypothesis~\ref{hyp:crucial} tells us that the overlap will be well-behaved: it
will break up nicely into disjoint copies of $\nbody{O}$ and $\nbody{W}$.}

In \cite{RSV07}, $\nbody{W}$ is always empty, so this hypothesis is not necessary.
\end{remark}

The final construction required for our main results is a map on symmetric
functions that will give an algebraic interpretation of the diagrammatic
operation $\circ_W$. 
We note that the definition below is the natural generalisation of \cite[Definition~7.18]{RSV07}.

\begin{definition}\label{def:importantmap}
Let $E$ and $W$ be skew diagrams such that $E=WOW$. Consider the\nnew{ }map of sets\no
$$
\begin{matrix}
\Lambda & \overset{(-) \circ_{W} s_E}{\longrightarrow}& \Lambda\no \\
    0  &      \longmapsto                               & 0\no \\
    f   &      \longmapsto                               & f \circ_W s_E
\end{matrix}
$$that consists of the composition of the following\no two maps $\Lambda \rightarrow \Lambda[t] \rightarrow \Lambda$
if $f\neq 0$.
If we think of $\Lambda$ as the polynomial algebra $\mathbb{Z}[h_1,h_2,\ldots]$, then we can
temporarily grade
$\Lambda$ and $\Lambda[t]$ by setting $\deg(t)=\deg(h_r)=1$ for all $r$. 
The first map $\Lambda \rightarrow \Lambda[t]$ then homogenises a polynomial in the $h_r$
with respect to the above grading, using the variable $t$ as the homogenisation variable.  

Meanwhile the second map is given by
$$
\begin{matrix}
\Lambda[t]& \longrightarrow &\Lambda\no \\
     h_r  & \longmapsto     &s_{E^{\amalg_W r}}\no \\
      t   & \longmapsto     &s_W.
\end{matrix}
$$
\end{definition}

For example, if $f=h_1h_2 h_3 - (h_3)^2 - h_2 h_4 + h_6$, then its image under the first
map is $h_1 h_2 h_3 -(h_3)^2 t - h_2 h_4 t + h_6 t^2$.  Therefore, 
\[
f \circ_W s_E = s_E s_{E^{\amalg_W 2}}s_{E^{\amalg_W 3}}
- (s_{E^{\amalg_W 3}})^2 s_W  - s_{E^{\amalg_W 2}}s_{E^{\amalg_W 4}} s_W
+s_{E^{\amalg_W 6}} (s_W)^2 .
\]

If $f = s_D$ for some skew diagram $D$, then we see that there is a nice way to think of $f \circ_W s_E$ in
terms of the Jacobi-Trudi decomposition matrix for $s_D$.  Specifically, we homogenise
by writing each entry of the form $s_\emptyset$ in the Jacobi-Trudi decomposition matrix as $h_0$.  
Then we replace $s_{(r)}=h_r$ by $E^{\amalg_W r}$ for $r \geq 0$. 
With the convention that $E^{\amalg_W 0} = W$, we now have that
$s_D \circ_W s_E$ is simply the
determinant of the resulting matrix.  The reader is invited to check that the example
above corresponds to this rule applied to the case of $f=s_D$ with $D=(4,2,2)/(1,1)$.
Consequently, we have

\begin{lemma}\label{lem:importantmap}
Let $\phi_1, \ldots, \phi_r$ be the ribbons of the Jacobi-Trudi decomposition of a 
skew diagram $D$.  For skew diagrams $E$ and $W$ with $E=WOW$, we have
\[
s_D \circ_W s_E
=  \left( \det \left(s_{(\phi_i \# \phi_j)}\right)_{i,j=1}^r \circ_W s_E \right) 
= \det \left(s_{(\phi_i \# \phi_j)\gencomp_W E}\right)_{i,j=1}^r, 
\]
where $\emptyset \gencomp_W E$ is defined to be $W$, and where 
$(\no\mathit{undefined} \gencomp_W E)$\no
is undefined.
\end{lemma}

Let $\upbody{D}$ denote the subset of elements of $D$ that have another element of $D$ one position
to their south.  Notice that $|\upbody{D}| = |\nwbody{D}|+r-1$, where $r$ is the number of rows in 
$D$.  
For symmetric functions $f$ and $g$, we will write $f=\pm g$ to mean that either
$f=g$ or $f=-g$.  

We are finally ready to put everything together and start reaping the 
rewards of our hard work.

\begin{conjecture}\label{con:gencomp}
For any skew diagram $D$, and a skew diagram $E$ satisfying 
Hypotheses~\ref{ite:maximal} -- \ref{ite:crucial}, we have
\begin{equation}\label{equ:gencompconj}
s_{D \gencomp_W E} \left(s_{\nbody{W}}\right)^{|\upbody{D}|} \left(s_{\nbody{O}}\right)^{|\nwbody{D}|}
= \pm \left( s_D \circ_W s_E \right).
\end{equation}
The sign on the right-hand side is a plus sign if $E=\rightright{W}{O}$ or 
$E=\upup{W}{O}$, and otherwise depends only on $D$.
Furthermore, if $E$ does not satisfy Hypothesis~\ref{ite:crucial}, 
then there exists a skew diagram $D$ for
which \eqref{equ:gencompconj} fails to hold.
\end{conjecture}

We can prove \eqref{equ:gencompconj} when $E$ is of the form 
$\rightup{W}{O}$ or $\upright{W}{O}$.  However, our proof techniques require one
further assumption for the two other forms of $E$.

\begin{hypothesis}\label{hyp:final}
If $E=\rightright{W}{O}$ or $E=\upup{W}{O}$,\ok then we assume that:
\begin{enumerate}
\addtocounter{enumi}{4}
\renewcommand{\theenumi}{\Roman{enumi}}
\item \label{ite:final} In $E$, at least one copy of $W$ has just one cell adjacent to $O$.
\end{enumerate}
\end{hypothesis}

\begin{theorem}\label{thm:gencomp}
For any skew diagram $D$,\no and a skew diagram $E$ satisfying 
Hypotheses~\ref{ite:maximal}  -- \ref{ite:final}, we have
\begin{equation}\label{equ:gencomp}
s_{D \gencomp_W E} \left(s_{\nbody{W}}\right)^{|\upbody{D}|} \left(s_{\nbody{O}}\right)^{|\nwbody{D}|}
= \pm\left( s_D \circ_W s_E \right).
\end{equation}
The sign on the right-hand side is a plus sign if $E=\rightright{W}{O}$ or 
$E=\upup{W}{O}$, and otherwise depends only on $D$.
\end{theorem}

Due to its length, we postpone the proof and devote the next section to it. 

\begin{remark}
Again, we can compare what appears here with the relevant parts of \cite{RSV07}.
If $W$ is empty and either $D$ or $E$ is a ribbon, we obtain
\begin{equation}\label{equ:rsv1}
s_{D \gencomp_\emptyset E} = s_D \gencomp_\emptyset s_E,
\end{equation}
which is their (7.2) and Proposition 7.5.
If $D$ and $W$ are ribbons with $\nbody{W}=\emptyset$, 
and $E=\rightright{W}{O}$ or $E=\upup{W}{O}$, we obtain 
\begin{equation}\label{equ:rsv2}
s_{D \gencomp_W E} = s_D \gencomp_W s_E,
\end{equation}
which is equivalent to \cite[Theorem~7.20]{RSV07}.
In the case that $E$ is a ribbon of the form 
$E=\rightup{W}{O}$ or $E=\upright{W}{O}$, we obtain a result
that has no analogue in \cite{RSV07}, but which implies the ribbon
staircase equivalence in its\ok Theorem 7.30. 
We will say more about this in Remark~\ref{rem:staircase}.
\end{remark}

\begin{remark}
We can think of 
$\left(s_{\nbody{W}}\right)^{|\upbody{D}|} \left(s_{\nbody{O}}\right)^{|\nwbody{D}|}$
as the term introduced by the overlaps in $D \gencomp_W E$.
This is consistent with \eqref{equ:rsv1} and \eqref{equ:rsv2}, and with
our observation that $\nbody{W}$ and $\nbody{O}$ are the 
shapes that form the overlap of {\boldmath{$E$}} and a copy of 
{\boldmath{$E$}} translated one position northwest.
\end{remark}

While the major proof in this article is that of Theorem~\ref{thm:gencomp}, 
our main target has been the following result, which serves as a mechanism for building
skew-equivalences.

\begin{theorem}\label{thm:equivs}
Suppose we have skew diagrams $D$, $D'$ with $D \sim D'$, and $E=WOW$ satisfying 
Hypotheses~\ref{ite:maximal}  -- \ref{ite:final}.  Then
\begin{equation}\label{equ:gencompcor}
D' \gencomp_W E \ \sim\  D \gencomp_W E \ \sim\  D \gencomp_{W^*} E^*.
\end{equation}
\end{theorem}

\begin{proof}
If $D \sim D'$, then $s_D = s_{D'}$ and so $s_D \circ_W s_E = s_{D'} \circ_W s_E$.
Since the left-hand side of \eqref{equ:gencomp} is Schur-positive, the sign 
on the right-hand side of \eqref{equ:gencomp} is determined by $s_D \gencomp_W s_E$.
By Theorem~\ref{thm:necconds},  
$|\nwbody{D}| = |\nwbody{D'}|$ and $|\upbody{D}| = |\upbody{D'}|$.
Applying Theorem~\ref{thm:gencomp}
then yields the first equivalence.

The second equivalence follows from the first and from Proposition~\ref{pro:rotation} and 
Lemma~\ref{lem:gencomp}(\ref{ite:star}) since
\[
D \gencomp_W E \,\sim\, (D \gencomp_W E)^* = D^* \gencomp_{W^*} E^*\,\sim\, D \gencomp_{W^*} E^*.
\]
\end{proof}

\begin{example}
\cite[Section~9]{RSV07} contains a list of the 6 skew-equivalences involving skew diagrams with at
most 18 cells that are not explained by the results there.  Using the equivalences of
Theorem~\ref{thm:equivs}, we can now explain these equivalences.  
In all cases, let $D=\begin{smallmatrix} \ti & \ti \\ \ti \end{smallmatrix}$ and $D'=D^*$.
Letting
\[
E = \begin{matrix}  & w & w  \\ \ti & \ti  \\ w & w \end{matrix}
\]
the first equivalence of \eqref{equ:gencompcor} gives
\[
\begin{matrix}
&&&& \ti & \ti \\
&& \ti & \ti & \ti \\
& \ti & \ti & \ti & \ti \\
& \ti & \ti \\
\ti & \ti \\
\ti & \ti 
\end{matrix}
\sim
\begin{matrix}
&&&& \ti & \ti \\
&&& \ti & \ti \\
&&& \ti & \ti \\
& \ti & \ti & \ti \\
\ti & \ti & \ti & \ti \\
\ti & \ti 
\end{matrix},
\]
which is the first of the 6 equivalences.
With 
\[
E = \begin{array}{cccc} \ti & \ti & w & w \\ w & w \end{array}
\]
we get the second equivalence, which is
\[
\begin{array}{cccccccc}
&&& \ti & \ti & \ti & \ti & \ti \\
&& \ti & \ti & \ti & \ti \\
\ti & \ti & \ti & \ti \\
\ti & \ti 
\end{array}
\sim
\begin{array}{cccccccc}
&&&& \ti & \ti & \ti & \ti \\
& \ti & \ti & \ti & \ti & \ti \\
\ti & \ti & \ti & \ti \\
\ti & \ti 
\end{array}.
\]
The remaining equivalences are obtained by setting $E$ to be
\[
\begin{matrix}  && w & w  \\ \ti & \ti & \ti \\ w & w \end{matrix}\ ,\ \ \ 
\begin{matrix}  & \ti & w & w \\ \ti & \ti \\ w & w \end{matrix}\ , \ \ \ 
\begin{matrix}  &&&  w  \\ w & \ti & \ti & w \\ w & \ti \end{matrix}\ \mbox{\ \ and\ \ } \ \ 
\begin{matrix} \ti & \ti & \ti & w & w \\ w & w \end{matrix}\ok
\]
respectively.
\end{example}

\begin{remark}\label{rem:staircase}
The case when $E$ is a ribbon of the form $\rightup{W}{O}$ or $\upright{W}{O}$
falls under the ribbon staircase construction of \cite{RSV07}.
Theorem 7.30 there amounts to the equivalence 
$D^* \gencomp_W E \sim D \gencomp_W E$.  Therefore, the first equivalence
of \eqref{equ:gencompcor}, even when $E$ is just a ribbon, shows that
more general ribbon staircase equivalences exist.  As a first example, one
could take $E=\begin{smallmatrix} \ti & w \\ w \end{smallmatrix}$ and 
$D \sim D'$ to be the unique non-trivial skew-equivalence for skew diagrams
with at most 8 cells, i.e.
\[
E \gencomp_{(1)} E \ \ =\ \  
\begin{matrix} && \ti & \ti \\ \ti & \ti & \ti \\ \ti & \ti \\ \ti \end{matrix}
\ \ \sim\ \ 
\begin{matrix} & \ti & \ti & \ti \\ & \ti & \ti \\ \ti & \ti \\ \ti \end{matrix}
\ \ =\ \  E^* \gencomp_{(1)} E.
\]
We leave the construction of $D \gencomp_W E$ and $D' \gencomp_W E$, each 
of which has 20 cells, as an exercise.
\end{remark}

\section{Proof of Theorem~\ref{thm:gencomp}}
\label{sec:proof}
We devote this section to our main proof.  As we will see, the proof
divides into two cases: $E=\rightright{W}{O}$ and
$E=\rightup{W}{O}$.  The first thing we will do is show that proofs of these two
cases respectively\no imply the result for $E=\upup{W}{O}$ and $E=\upright{W}{O}$.\ok  The proof of each case has two major steps.  The first is to define 
an appropriate outside ribbon decomposition of $D \gencomp_W E$.  The second
step is to apply Sylvester's Determinantal Identity to the resulting Hamel-Goulden
matrix\ok and to\ok show that \eqref{equ:gencomp} results.  

\subsection{Reduction to two cases}
\nnew{If $W=\emptyset$, then we are in the case $E=\rightright{W}{O}$, which is one of the cases we will prove explicitly; therefore, assume that $W \neq \emptyset$.}
\new{For the following reasons,}
we see that \eqref{equ:gencomp} holds for $D$ and $E$ if and only if the 
corresponding equation for $D^*$ and $E^*$ holds:
\begin{itemize}
\item By Proposition~\ref{pro:rotation} and 
Lemma~\ref{lem:gencomp}(\ref{ite:star}), $s_{D \gencomp_W E} = s_{D^* \gencomp_{W^*} E^*}$. 
\item We have $E^* = W^*O^*W^*$, and we see that
$\nbody{O^*} = \left(\nbody{O}\right)^*$ and $\nbody{W^*} = \left(\nbody{W}\right)^*$.
\item $|\upbody{D}| =  |\upbody{D^*}|$ and $|\nwbody{D}| = |\nwbody{D^*}|$.  
\item By Proposition~\ref{pro:rotation} and since $\left(E^{\amalg_W n}\right)^* = (E^*)^{\amalg_{W^*} n}$, $s_D \circ_W s_E = s_{D^*} \circ_{W^*} s_{E^*}$.
\end{itemize}
\new{Moreover, applying the $\omega$ involution 
that sends $s_D$ to $s_{D^t}$,
we get that \eqref{equ:gencomp} holds for $D$ and $E$ if and only if
\begin{equation}\label{equ:afteromega}
s_{(D \gencomp_W E)^t} \left(s_{\left(\nbody{W}\right)^t}\right)^{|\upbody{D}|} \left(s_{\left(\nbody{O}\right)^t}\right)^{|\nwbody{D}|}
= \pm \left( s_D \circ_{W^t} s_{E^t}\right)  .
\end{equation}
Applying a similar argument to that for $D^*$ and $E^*$, we get that \eqref{equ:afteromega}
holds if and only if
\[
s_{D^* \gencomp_{W^t} E^t} \left(s_{\nbody{W^t}}\right)^{|\upbody{D^*}|} \left(s_{\nbody{O^t}}\right)^{|\nwbody{D^*}|}
= \pm \left( s_{D^*} \circ_{W^t} s_{E^t} \right) ,
\]
i.e.\ \eqref{equ:gencomp} holds for $D$ and $E$ if and only if the 
corresponding equation for $D^*$ and $E^t$ holds.  }
Therefore, \eqref{equ:gencomp} holds for $E$ and all skew diagrams $D$ if and only 
if it also holds with $E^*$, $E^t$ or $(E^*)^t$ in place of $E$, and all skew diagrams $D$.  
Furthermore, $E^*$, $E^t$ and $(E^*)^t$ satisfy Hypotheses~\ref{ite:maximal}  -- \ref{ite:final}
if and only if $E$ does.
Therefore, we can assume that $E$ is arranged in a convenient manner:
\begin{enumerate}
\item \label{ite:amalgcomp}
If $E=\rightright{W}{O}$ or $E=\upup{W}{O}$,\ok then let us choose $E$
so that $E = \rightright{W}{O}$ and $\Wsw$ has at least as many rows
in common with $O$ as $\Wne$.    By Hypothesis~\ref{ite:final}, $\Wne$
thus has just one row in common with $O$. 
\item \label{ite:staircase}
 If $E=\rightup{W}{O}$ or $E=\upright{W}{O}$, we will assume that $E=\rightup{W}{O}$.
\end{enumerate}
The advantage of these forms of $E$
is that every copy of $\nbody{W}$ is contained in the corresponding copy of $W$.

Let us also reduce the second assertion in Theorem~\ref{thm:gencomp} to
these two cases. 
Since the left-hand side of \eqref{equ:gencomp} is Schur-positive, the sign
on the right-hand side depends only on $s_D \gencomp_W s_E$.
Therefore, it follows from our argument above that
the signs on the right-hand side in front of 
$s_D \gencomp_W s_E$ and
$s_{D^*}  \gencomp_{W^t} s_{E^t}$  will
be the same.  We know that if $E=\upup{W}{O}$, then $E^t = \rightright{W^t}{O^t}$.
Therefore, we get a plus sign on the right-hand side for every skew diagram $D$
and every $E=\upup{W}{O}$\no satisfying Hypotheses~\ref{ite:maximal}  -- \ref{ite:final},
if and only if the same applies when $E=\rightright{W}{O}$\no satisfying the hypotheses.
It also follows from our argument above that the signs on the right-hand side in
front of $s_D \gencomp_W s_E$ and $s_D \gencomp_{(W^*)^t} s_{(E^*)^t}$ will be the same.
If $E=\upright{W}{O}$,\ok then $(E^*)^t = \rightup{(W^*)^t}{(O^*)^t}$.  Hence, for fixed $D$, if the 
sign on the right-hand side of \eqref{equ:gencomp} is the same for all 
$E$ of the form $\rightup{W}{O}$ satisfying the hypotheses, then it will be the same for all $E$ of the
forms $\rightup{W}{O}$ and $\upright{W}{O}$ satisfying the hypotheses.

It will now be necessary to treat Cases~\ref{ite:amalgcomp} and \ref{ite:staircase} above separately.

\subsection{The outside decomposition for Case~\ref{ite:amalgcomp}}
Our first major task is to construct an appropriate ribbon decomposition for $D \gencomp_W E$. For $d \in D$ and $e \in E$, 
let us write $d \gencomp_W e$ to denote the cell that corresponds to $e$ in the $E_d$ copy
of $E$ in $D \gencomp_W E$.   Then until further notice, 
for $D' \subseteq D$ and $E' \subseteq E$, we define
$D' \gencomp_W E'$ to be the subset of $D \gencomp_W E$ consisting of those cells
of the form $d \gencomp_W e$, where $d \in D'$ and $e \in E'$.  
We construct a ribbon decomposition
for $D \gencomp_W E$ by first constructing a ribbon decomposition for $D$.  
As mentioned in Remark~\ref{rem:horizdecomp}, we will 
take the Jacobi-Trudi decomposition so that our ribbons are simply the rows
$\phi_1, \ldots, \phi_r$ from north to south of $D$.  Let $\Phi$ denote
the cutting strip of this decomposition, which is just a single row of size $|\nwD|$.
A ribbon decomposition for
$\phi_1 \gencomp_W E$ is now suggested by Definition~\ref{def:nbody}.
Indeed, every copy of $E$ can be partitioned into two copies of $\nbody{W}$, one copy
of $\nbody{O}$ and a ribbon $\south{E}$ which runs along the south border of $E$.
Extending this idea, $\phi_1 \gencomp_W E$ can be
partitioned into $|\phi_1|+1$ copies of $\nbody{W}$, $|\phi_1|$ copies of $\nbody{O}$ and
a ribbon $\south{\phi_1 \gencomp_W E}$ obtained by amalgamating the $|\phi_1|$ 
copies of $\south{E}$.  
Construct the southeast decompositions for the copies of $\nbody{W}$ and $\nbody{O}$ to 
obtain a ribbon decomposition for $\phi_1 \gencomp_W E$.
For example, if 
\[
D = \begin{array}{cc} \ti & \ti \\ \ti & \ti \\ \ti \end{array} \mbox{\ \ and\ \ } 
E = \begin{array}{cccccc} 
&&&& w & w \\
&& \ti & \ti & w & w \\
w & w & \ti & \ti \\
w & w & \ti
\end{array},
\]
then $\phi_1 \gencomp_W E$ has a ribbon decomposition given by the ribbons labelled 
$1$, $8$, $9$, $10$, $11$, $12$ in Figure~\ref{fig:Case1all}.
\begin{figure}\[
D \gencomp_W E = 
\begin{array}{cccccccccccccccc}
&&&&&&&&&&&&&& 12 & 12 \\
&&&&&&&&&&&& 11 & 1 & 1 & 1 \\
&&&&&&&&&& 10 & 10 & 1 & 1 \\
&&&&&&&& 9 & 1 & 1 & 1 & 1\\
&&&&&& 8 & 8 & 1 & 1 & 2 & 2 & 2 \\
&&&&&& 1 & 1 & 1 & 2 & 2 \\
&&&&& 7 & 2 & 2 & 2 & 2 \\
&&& 6 & 6 & 2 & 2 \\
&&& 2 & 2 & 2 \\
&& 5 & 3 & 3 & 3 \\
4 & 4 & 3 & 3 \\
3 & 3 & 3
\end{array}
\]\caption[]{}\label{fig:Case1all}\end{figure}

Now consider what happens when we add in the next part of $D \gencomp_W E$, namely
$\phi_2 \gencomp_W E$.  Some of $\phi_2 \gencomp_W E$ will have already been included as
elements of $\phi_1 \gencomp_W E$, while some will be new. 
If $d'$ is one position northwest of $d$ in $D$, then it follows
from Definition~\ref{def:gencomp} that $E_{d'}$ will 
be one position northwest of $E_d$ in $D \gencomp_W E$.  Furthermore, we see that
the $\nbody{O}$ and the two copies of\ok $\nbody{W}$ are exactly the subset of $E_d$ that
will be contained in $E_{d'}$.
In particular, this implies that the ribbon
$\phi_2 \gencomp_W \south{E}$ will be new.  Also, for $d \in \phi_2$, 
the cells of the 
copy of $\nbody{O}$ in $d \gencomp_W E$ will be new if and only if
there is no cell $d'$ one position northwest of $d$ in $D$,\ok i.e.\ if and only if $d \in \nwD$.
The same applies to the cells of the lower copy of $\nbody{W}$ in $d \gencomp_W E$.  
The cells of the upper copy of $\nbody{W}$ in $d \gencomp_W E$ will be new if and only
if there is no cell $d'$ one position north of $d$ in $D$.  We conclude that $\phi_2 \gencomp_W E$
contributes $|\phi_2 \cap \nwD|$ new copies of both $\nbody{O}$ and $\nbody{W}$ to 
$D \gencomp_W E$.  

Continuing in this manner, $D \gencomp_W E$ can be partitioned into the ribbons
$\phi_i \gencomp_W \south{E}$ for $i=1,\ldots,r$, 
along with $|\nwD|$ copies of $\nbody{O}$ and $|\nwD|+1$ 
copies of $\nbody{W}$.  Due to Hypotheses~\ref{ite:separation} and \ref{ite:crucial}, 
the copies of $\nbody{O}$
and $\nbody{W}$ are all pairwise non-adjacent.\no
Using Hypothesis~\ref{ite:crucial}, 
one can also check that these copies of $\nbody{O}$ and $\nbody{W}$, as well as the
endpoints of the ribbons $\phi_i \gencomp_W \south{E}$, are on the outside of $D \gencomp_W E$
in the appropriate sense.
Therefore, we can construct an outside ribbon 
decomposition $\theta_1, \ldots, \theta_N$ for $D \gencomp_W E$ as follows:
\begin{itemize}
\item For $i=1,\ldots,r$, let $\theta_i := \phi_i \gencomp_W \south{E}$.
\item Construct the southeast decomposition for each of the $|\nwD|+1$ copies of $\nbody{W}$
and each of the $|\nwD|$ copies of $\nbody{O}$.  Label the resulting ribbons 
$\theta_{r+1},\ldots, \theta_N$ from 
southwest to northeast according to their southwest endpoints.  
\end{itemize}
We let $\Theta$ denote the cutting strip of this ribbon decomposition.
Figure~\ref{fig:Case1all}  is an example of such a decomposition.

\subsection{Applying Sylvester's Determinantal Identity in Case~\ref{ite:amalgcomp}}
Now that our outside ribbon decomposition of $D \gencomp_W E$ is defined, we 
proceed to the second part of the proof for $E=\rightright{W}{O}$.
We wish to apply Sylvester's Determinantal
Identity to $M$, the Hamel-Goulden matrix for our decomposition $\theta_1, \ldots, \theta_N$.
We let the $S$ in Theorem~\ref{thm:sylvester} be the set
$\{r+1, \ldots, N\}$.  Then $M[S,S]$ is block lower-triangular,\ok and the blocks on the diagonal
alternate between Hamel-Goulden matrices for $\nbody{W}$ and Hamel-Goulden matrices
for $\nbody{O}$.
The left-hand side of \eqref{equ:sylvester} is thus
\begin{equation}\label{equ:lhs}
s_{D \gencomp_W E} \left(\left(s_{\nbody{W}}\right)^{|\nwD|+1}\left(s_{\nbody{O}}\right)^{|\nwD|}\right)^{r-1}.
\end{equation}

We next wish to evaluate $M[S\cup\{i\}, S\cup\{j\}]$ for $i,j \in \{1,2,\ldots,r\}$.  
The top-left entry
of this submatrix is $s_{\theta_i \# \theta_j}$.  
The remaining entries of the first row are $s_{\theta_i \# \theta_l}$
for $l =r+1, \ldots, N$.  
The remaining entries of the first column are 
$s_{\theta_k \# \theta_j}$
for $k =r+1, \ldots, N$.  The $(k,l)$-entry for $k,l>1$ is $s_{\theta_{k+r-1} \# \theta_{l+r-1}}$.
Therefore, our goal of evaluating $M[S\cup\{i\}, S\cup\{j\}]$ has the following two parts:
\begin{enumerate}
\renewcommand{\theenumi}{\alph{enumi}}
\item \label{ite:defined} If $\theta_i \# \theta_j$ is defined, then 
$s_{\theta_i \# \theta_l} = s_{(\theta_i \# \theta_j) \# \theta_l}$ and
$s_{\theta_k \# \theta_j} = s_{\theta_k \# (\theta_i \# \theta_j)}$.  Therefore, 
we wish to show that a skew diagram 
$F$ with the following property exists:
when using the cutting strip $\Theta$, 
the ribbons of the resulting decomposition of $F$ correspond to the same portions of  
$\Theta$ as $$\theta_i \# \theta_j, \theta_{r+1}, \theta_{r+2}, \ldots, \theta_N.$$
Then $M[S\cup\{i\}, S\cup\{j\}]$ for $i,j \in \{1,2,\ldots,r\}$ will be the Hamel-Goulden matrix for $F$.
\item \label{ite:null} If $\theta_i \# \theta_j$ is undefined, then we wish to show that
$\det M[S\cup\{i\}, S\cup\{j\}] = 0$.
\end{enumerate}

We first show (\ref{ite:null}).  
{In our running example, we can take $i=3$ and $j=1$. }
We will show that, for some $t$, the submatrix of $M[S\cup\{i\}, S\cup\{j\}]$ 
consisting of its first $t+1$ rows has rank less than $t+1$, implying
the result.  
As usual, let $p(\theta_k)$ and $q(\theta_k)$ respectively denote the starting contents and
ending contents of the ribbon $\theta_k$. 
We know that $\theta_k \# \theta_l$ is
defined if and only if $p(\theta_l) \leq q(\theta_k) + 1$. 
We also know that $q(\theta_i)$ is the content of the 
northeast endpoint $q_i$ of the ribbon $\phi_i \gencomp_W \south{E}$.
In particular, $q_i$ is in a copy $W_i$ of $W$.  
We will be comparing the positions of copies of $W$, $\nbody{W}$ and $\nbody{O}$ 
according to how \new{far} southwest or northeast they are.  We will say that a subdiagram $A$ of
$D \gencomp_W E$ is weakly southwest of a subdiagram $B$ if
the maximum content of the cells of $A$ is weakly less than the maximum
content of the cells of $B$.  As a variation of this, we will say that $A$ is
strictly northeast of $B$ if the minimum content in $A$ is strictly greater than the
maximum content in $B$.  In effect, we are comparing subdiagrams according to their positions
on the cutting strip of $D \gencomp_W E$.  
Let $t$ be the maximum positive integer such that $\theta_{r+t}$ is in the 
southeast decomposition of a copy of $\nbody{W}$ or $\nbody{O}$ that
is weakly southwest of $W_i$, or set $t=0$ if no such positive integer
exists.  
{In our running example where $i=3$, we have that $r+t=6$.}
{Due to} the way we have labelled
$\theta_{r+1}, \ldots, \theta_N$, this \new{definition of $t$} tells us that 
$\theta_{r+1}, \ldots, \theta_{r+t}$ are contained in copies of $\nbody{W}$
and $\nbody{O}$ weakly southwest of $W_i$, while
$\theta_{r+t+1}, \ldots, \theta_N$ are contained in copies strictly northeast
of $W_i$.
We have that:
\begin{itemize}
\item $p(\theta_j) > q(\theta_i)+1$, since $\theta_i \# \theta_j$ is undefined.
\item $p(\theta_l) > q(\theta_i)+1$ for $l=r+t+1, \ldots, N$.
\new{Indeed, suppose $p(\theta_l) \leq q(\theta_i)+1$, and let
$p_l$ denote the southwest endpoint of $\theta_l$.  By Hypothesis~\ref{ite:separation},
the only danger is if $\theta_l$ is contained in the copy $O_l$ of $O$ that
is immediately northeast of $W_i$. 
Since $\theta_i \# \theta_j$ is undefined, we know that $i \neq 1$, and we also know by definition
of $q_i$ 
that $W_i$ is the northeasternmost copy of $W$ in $\phi_i \gencomp_W E$.   Therefore, 
there exists a copy $W_l$ of $W$ that is a translation of $W_i$ one position to the northwest.
Observe that $W_l$ and $O_l$ are from the same copy of $E$, and 
let $q_l \in W_l$ denote the translation of $q_i$ one position to the northwest.
Since $q_i$ has no cell immediately to its east, we know that $q_l$ has a cell $o$ of $O_l$ one position to its east.  In fact, we know that $o \not\in \nbody{O_l}$,\ok so, in particular, $o \neq p_l$.
Since $p(\theta_l) \leq q(\theta_i)+1$,
$p_l$ must be strictly west and weakly north of $o$.  However, since $o, p_l \in O_l$ and 
$O_l$ is a skew diagram, this implies that $q_l$ is contained in 
$O_l$, a contradiction.} 

\item $p(\theta_j) > q(\theta_k)+1$, for $k=r+1, \ldots, r+t$.  Indeed, we know that
$p(\theta_j)$ is the content of the southwest endpoint $p_j$ of the ribbon
$\phi_j \gencomp_W \south{E}$.  Therefore, $p_j$ is an element
of a copy $W_j$ of $W$.  However, since 
$p(\theta_j) > q(\theta_i)+1$, $W_j$ must be strictly northeast of $W_i$, while
$\theta_k$ is weakly southwest of $W_i$.  Hypothesis~\ref{ite:separation}
then implies the claim.
\item $p(\theta_l) > q(\theta_k)+1$, for $l=r+t+1, \ldots, N$ and
$k=r+1, \ldots, r+t$.  Indeed, $\theta_l$ is contained in 
a copy of $\nbody{W}$ or $\nbody{O}$ that is strictly northeast of
{$W_i$}, while $\theta_k$ is contained in a copy of $\nbody{W}$ or
$\nbody{O}$ that is weakly southwest of $W_i$.  By Hypotheses~\ref{ite:separation}
and \ref{ite:crucial}, no two copies of $\nbody{W}$ or $\nbody{O}$ are
adjacent, implying the claim.
\end{itemize}
Therefore,  $\theta_k \# \theta_l$ is undefined, and so
$s_{\theta_k \# \theta_l}=0$, whenever
\new{$k \in \{i,r+1, r+2, \ldots, r+t\}$ and $l \in \{j, r+t+1, r+t+2, \ldots, N\}$.
\nnew{This shows that the submatrix of $M[S\cup\{i\}, S\cup\{j\}]$ consisting of its first $t+1$ rows has rank at most
$t$, as required for (b).}}

\new{To construct the diagram $F$ needed to establish (a), we move from the setting of $D \gencomp_W E$
to the setting of $\Phi \gencomp_W E$, recalling that $\Phi$ is the cutting strip for the
Jacobi-Trudi decomposition of $D$, i.e.\ $\Phi$ is a single row of size $|\nwD|$.}
As a result, 
$d \gencomp_W e$ now denotes the cell that corresponds to $e$ in the $E_d$ copy
of $E$ in $\Phi \gencomp_W E$.  Subsequently, 
a term of the form $\Phi' \gencomp_W E'$
with $\Phi' \subseteq \Phi$ and $E' \subseteq E$ denotes those 
$d \gencomp_W e \in \Phi \gencomp_W E$ such that $d \in \Phi'$ and $e \in E'$.
We will continue to use $\Theta$ as our cutting strip.  
Notice that this change does not affect the definition 
$\theta_i := \phi_i \gencomp_W \south{E}$, since we get the same portion of $\Theta$
in the setting of $\Phi \gencomp_W E$ as we did in the $D \gencomp_W E$ setting.
We have that
\[  
\theta_i \# \theta_j = (\phi_i \gencomp_W \south{E}) \# (\phi_j \gencomp_W \south{E})
= (\phi_i \# \phi_j) \gencomp_W \south{E}.
\]
If $\phi_i \# \phi_j = \emptyset$, then 
$(\phi_i \# \phi_j) \gencomp_W \south{E}$ should be defined to be the portion
of $W$ that is contained in $\south{E}$, namely $W \setminus \nbody{W}$.

As usual, let $\Phi[p(\phi_j), q(\phi_i)]$ denote the portion of the 
cutting strip $\Phi$ for $D$ corresponding to $\phi_i \# \phi_j$.  We will let the southwesternmost
cell of the cutting strip have content 1, so the cutting strip itself can
be expressed as $\Phi[1,|\nwD|]$.
Let $\nbody{E}$ denote $E \setminus \south{E}$, which we know in our case consists of
two copies of $\nbody{W}$ and a copy of $\nbody{O}$.  We write $\nbody{W}\nbody{O}$
(resp.\ $\nbody{O}\nbody{W}$)
to denote the result of deleting the upper (resp.\ lower) copy of $\nbody{W}$ from $\nbody{E}$.
When $\Phi[p(\phi_j), q(\phi_i)]$ is defined,  
we claim that the required skew diagram $F$ is the subset of $\Phi \gencomp_W E$ consisting of
\begin{equation}\label{equ:byconstruction}
\left(\Phi[p(\phi_j), q(\phi_i)] \gencomp_W \south{E}\right) \cup \left(\Phi \gencomp_W \nbody{E}\right) .
\end{equation}

In our running example, we could take $i=2$ and $j=3$.  
Then $F$ is the following skew diagram, where we label the cells according
to the label of the corresponding cell in the ribbon decomposition of $D \gencomp_W E$:\ok
\[
\begin{matrix}
&&&&&&&&&&&&&&&& 12 & 12 \\
&&&&&&&&&&&&&& 11 \\
&&&&&&&&&&&& 10 & 10 \\
&&&&&&&&&& 9 & 2 & 2 & 2 \\
&&&&&&&& 8 & 8 & 2 & 2 \\
&&&&&& 7 & 2 & 2 & 2 & 2 \\
&&&& 6 & 6 & 2 & 2 \\
&& 5 & 3 & 3,2 & 3,2 & 2 \\
4 & 4 & 3 & 3 \\
3 & 3 & 3 
\end{matrix}\nnew{.}
\]
Intuitively, we have deleted the ribbon in $D \gencomp_W E$ labelled 1, and then translated
some of the remaining ribbons to the southeast.
Since translations to the southeast do not affect the contents
of cells, for each remaining ribbon we are using the same portions of the cutting strip as before.
The union of \eqref{equ:byconstruction} is a disjoint union,\ok 
since $\Phi$ is a ribbon and $\south{E}$ is disjoint from $\nbody{E}$ in $E$.
It can be rewritten as
\begin{eqnarray*} 
& & 
\left(\Phi[p(\phi_j), q(\phi_i)] \gencomp_W \south{E}\right) 
\sqcup \left(\Phi[1,p(\phi_j)-1] \gencomp_W \nbody{W}\nbody{O}\right)  \\
& & \sqcup \left(\Phi[p(\phi_j), q(\phi_i)] \gencomp_W \nbody{E}\right)
\sqcup \left(\Phi[q(\phi_i)+1,|\nwD|] \gencomp_W \nbody{O}\nbody{W}\right) \\
& = &
\left( \Phi[1,p(\phi_j)-1] \gencomp_W \nbody{W}\nbody{O}\right)
\sqcup  \left(\Phi[p(\phi_j), q(\phi_i)] \gencomp_W E\right)  \\
& & \sqcup \left(\Phi[q(\phi_i)+1,|\nwD|] \gencomp_W \nbody{O}\nbody{W}\right). 
\end{eqnarray*}
We make three observations about the latter disjoint union:
\begin{itemize}
\item The three terms are (not necessarily connected)
skew diagrams when considered individually.
\item Their union is obtained from $\Phi \gencomp_W E$
by removing an initial portion and a final portion of the south border ribbon
$\Phi \gencomp_W \south{E}$.  
\item Using Hypotheses~\ref{ite:crucial} and \ref{ite:final}, we see that at least
one empty diagonal separates each pair of terms in the disjoint union.  
\end{itemize}
Taken together, these observations tell us that the disjoint union is a skew diagram.
By construction (see \eqref{equ:byconstruction}), 
using the cutting strip $\Theta$, 
the ribbons of the resulting decomposition correspond to the same portions of  
$\Theta$ as $\theta_i \# \theta_j, \theta_{r+1}, \theta_{r+2}, \ldots, \theta_N$. 
Thus we have shown (\ref{ite:defined}).

Evaluating the corresponding skew Schur function, we finally deduce that
\[
\det M[S \cup \{i\}, S \cup \{j\}] = 
s_{(\phi_i \# \phi_j)\gencomp_W E} \left(s_{\nbody{W}}s_{\nbody{O}}\right)^{|\nwD|-q(\phi_i)+p(\phi_j)-1},
\]
since $s_{\nbody{O}\nbody{W}} = s_{\nbody{W}\nbody{O}} = s_{\nbody{W}}s_{\nbody{O}}$.

Therefore, by Lemma~\ref{lem:importantmap}, 
$\det \syl{M}{S}$ evaluates to

\new{
\[
\left(s_{\nbody{W}}s_{\nbody{O}}\right)^I \det \left(s_{(\phi_i \# \phi_j)\gencomp_W E}\right)_{i,j=1}^r 
=  \left(s_{\nbody{W}}s_{\nbody{O}}\right)^I \left( s_D \circ_W s_E \right), 
\]
}
where 
\[
I = r(|\nwD|-1) - \sum_{i=1}^r q(\phi_i) + \sum_{j=1}^r p(\phi_j).
\]
Plugging this and \eqref{equ:lhs} into \eqref{equ:sylvester}, we obtain
\[
s_{D \gencomp_W E} \left(s_{\nbody{W}}\right)^{(|\nwD|+1)(r-1)-I} \left(s_{\nbody{O}}\right)^{|\nwD|(r-1)-I} = s_D \circ_W s_E .
\]
Since
\[
|\nwD|(r-1)-I = - |\nwD| + \sum_{i=1}^r (q(\phi_i)-p(\phi_i)+1) = |D| - |\nwD| = |\nwbody{D}| 
\]
and 
\[
(|\nwD|+1)(r-1)-I = |\nwbody{D}| + r-1 = |\upbody{D}| ,
\]
we conclude that 
\[
s_{D \gencomp_W E} \left(s_{\nbody{W}}\right)^{|\upbody{D}|} \left(s_{\nbody{O}}\right)^{|\nwbody{D}|}
= s_D \circ_W s_E,
\]
as required.

\subsection{The outside decomposition for Case~\ref{ite:staircase}}
Throughout, it will be advantageous to think of $D \gencomp_E W$ in terms of 
Definition~\ref{def:gencomp2}.  
In particular, we know that when $d$ is one position south of $d'$ on the same ribbon in $D$,
then $E_d$ and $E_{d'}$ appear in the form $E_d \amalg_W E_{d'}$.\ok
It is necessary to explicitly 
add an extra copy of $W$ one position southeast of $E_d \cap E_{d'}$ if and only if
both $d$ and $d'$ are elements of $\se{D}$.
Accounting for these extra copies of $W$ is
what makes the notation of this case difficult.  In Case~\ref{ite:amalgcomp},
we wrote $d \gencomp_W e$ to denote the cell that corresponds to $e$ in the
$E_d$ copy of $E$ in $D \gencomp_W E$.  However, 
when $E=\rightup{W}{O}$, the
cells of the extra copies of $W$ do not naturally take the form $d \gencomp_W e$
for any $d$ and $e$.  We begin by concocting a way to rectify this situation.

Recall that $\emptyset \gencomp_W E = W$.  So we could think of these extra copies
of $W$ as being the contribution from some empty ribbons that we add to the
northwest decomposition of $D$.  Specifically, if $d$ is one position south of $d'$ on the
same ribbon in $D$ and $d, d' \in \se{D}$, then we can add imaginary cells 
$i(d)$ and $i(d')$ one position southeast of $d$ and $d'$ respectively.  Then to the
northwest decomposition we add an empty ribbon which starts at $i(d')$ and
ends at $i(d)$.  When we have added the empty ribbons for all such $d$ and $d'$,
we call the resulting ribbon decomposition the $\emph{enhanced northwest \new{decomposition}}$.
The big advantage of this ribbon decomposition is that 
every ribbon of size $l$ in the enhanced northwest decomposition contributes
$E^{\amalg_W l}$ to $D \gencomp_W E$, where we define $E^{\amalg_W 0}$ to be $W$\ok
and where\ok $D \gencomp_W E$ is exactly the union of these contributions, suitably positioned.
For $D$ in Example~\ref{exa:extraws}, we would have
\begin{equation}\label{equ:enhanced}
\begin{matrix}
& a_5 \\
a_3 & a_4 & i(a_5) \\
a_2 & b & i(a_4) \\
a_1 & i(a_2) \\
& i(a_1)
\end{matrix}\nnew{\ } .
\end{equation}
Let $p(\phi_i)$ and
$q(\phi_i)$ respectively denote the contents of the southwest and northeast endpoints of
a ribbon $\phi_i$.
We leave it as a nice exercise for the reader to check that the enhanced ribbon decomposition
of any $D$ has $r$ ribbons, where $r$ is the number of rows of $D$.  
Therefore, let $\phi_1, \ldots, \phi_r$ denote this set of ribbons, ordering the ribbons
so that $q(\phi_1) \geq q(\phi_2) \geq \cdots \geq q(\phi_r)$.  
In particular, we will have $\phi_1 = \nwD$.  

We are now in a good position to describe a ribbon decomposition for $D \gencomp_W E$.
First let us give an example to which the reader can refer for intuition.
Suppose that
\[
D =
\begin{matrix}
\ti & \ti \\ \ti & \ti \\ \ti 
\end{matrix}
\mbox{\ \ and\ \ }
E = \begin{matrix}
&&& w & w \\ &&& w \\ && \ti & \ti \\ w & w & \ti & \ti \\ w & \ti & \ti 
\end{matrix}.
\]
The ribbon decomposition of $D \gencomp_W E$ that we will construct is
\[
\begin{matrix}
&&&&&&&&&&&& 12 & 12 \\
&&&&&&&&&&&& 1 \\
&&&&&&&&&&& 11 & 1 \\
&&&&&&&&& 10 & 10 & 1 & 1 \\
&&&&&&&&& 1 & 1 & 1 \\
&&&&&&&& 9 & 1 & 2 \\
&&&&&& 8 & 8 & 1 & 1 & 2 \\
&&&&&& 1 & 1 & 1 & 2 & 2 \\
&&&&& 7 & 1 & 2 & 2 & 2 \\
&&& 6 & 6 & 1 & 1 \\
&&& 1 & 1 & 1 \\
&& 5 & 1 & 3 \\
4 & 4 & 1 & 1 \\
1 & 1 & 1
\end{matrix}.
\]

We must give a precise description of the general form of this ribbon decomposition.
We will divide $D \gencomp_W E$ into two subdiagrams, $L$ and $U$, each of
which is relatively easy to handle.  We will then give $L$ a northwest decomposition
and $U$ a decomposition that is close to a southeast decomposition.
These two decompositions will be compatible, resulting in a global outside
ribbon decomposition for $D \gencomp_W E$. 

The upper subdiagram $U$ is
the contribution from $\nwD$.  In other words, $U$ is the unique
subdiagram of shape $E^{\amalg |\nwD|}$ in $D \gencomp_W E$ that
includes $\nw{D \gencomp_W E}$.  In our example, $U$ is the
subdiagram consisting of those cells with labels from the set $\{1, 4, 5, 6, \ldots, 12\}$.

Now consider the contribution to $D \gencomp_W E$
from those ribbons of the \new{enhanced northwest
decomposition} of $D$, other than $\nwD$.
If $\phi_i$ has size $l$, then we noted that
$\phi_i$ contributes the skew diagram $E^{\amalg_W l}$ to $D \gencomp_W E$. 
However, much of this contribution will be contained in the contribution 
of $\phi_j$, where $\phi_j$ is immediately to the northwest of $\phi_i$.  
\nnew{As before, define $\south{E}$} as the result of removing the two copies of $\nbody{W}$ and the
copy of $\nbody{O}$ from $E$.  

For $i=1,\ldots,r$,\no in the skew diagram $E^{\amalg_W l}$ contributed by $\phi_i$, let
$\phi_i \gencomp_W \south{E}$ denote the subdiagram
obtained by considering only $\south{E}$ in each copy of $E$.
For example, 
$\phi_1 \gencomp_W \south{E} = \nwD \gencomp_W \south{E}$ is the
ribbon labelled 1 in our example.  It is clear that
$\phi_i \gencomp_W \south{E}$ is a ribbon for $i=1,\ldots,r$.
We can also see that
$\phi_1 \gencomp_W \south{E}, \ldots, \phi_r \gencomp_W \south{E}$
are disjoint in $D \gencomp_W E$, in part because $\phi_1, \ldots, \phi_r$ are
disjoint in $D$.
Furthermore, if $i>1$, then 
$\phi_i \gencomp_W \south{E}$ is exactly the contribution from $\phi_i$ that is not
contained in the contribution of $\phi_j$, where $\phi_j$ is immediately northwest
of $\phi_i$ in $D$.  
We let the lower subdiagram $L$ be
the union of $\phi_1 \gencomp_W \south{E}, \ldots, \phi_r \gencomp_W \south{E}$.  
One could think of $L$ as $D \gencomp_W \south{E}$.
In our example, $L$ is the subdiagram with cells\no labelled 1, 2, 3.

In our example, $U \cap L$ is the ribbon labelled 1.
In general, we see that  $U \cap L$ is the ribbon $\phi_1 \gencomp_W \south{E}$.
We also have that $U \cup L = D \gencomp_W E$.  Indeed, 
the empty ribbons from the
enhanced northwest decomposition contribute the cells that come from
the extra copies of $W$.  If $\phi_i$ is empty, by definition
$\phi_i \gencomp_W \south{E}$ is going to be the set of cells of $W$ that are
in $\south{E}$, i.e.\ $W \setminus \nbody{W}$.  We can now see that 
if a cell of $D \gencomp_W E$
is not in $\phi_i \gencomp_W \south{E}$ for some $i$, then it must be an
element of $E_d$ for some $d \in \nwD$, and hence is in $U$.

We can now define the ribbon decomposition 
$\theta_1, \ldots, \theta_N$ of $D \gencomp_W E$.  
For $i=1,\ldots,r$, we let $\theta_i = \phi_i \gencomp_W \south{E}$.  
Then, by definition of $L$, $\theta_1, \ldots, \theta_r$ give a ribbon decomposition
of $L$.  Certainly, if $\phi_i$ is empty, then $\theta_i$ is an outside ribbon in
$L$. 
If $\phi_i$ is non-empty,\no 
then $\theta_i$ will be an outside ribbon in $L$ because
$\phi_i$ was an outside ribbon in $D$.  Furthermore, since $\theta_1$ is $\nw{L}$,
it must be the cutting strip of this outside decomposition of $L$, and so
$\theta_1, \ldots, \theta_r$ give the northwest ribbon decomposition of $L$.  
One can also check that $\theta_1, \ldots, \theta_r$ will be outside ribbons in 
$D \gencomp_W E$.

We proceed to decompose $U$ into ribbons in such a way that we get
a global decomposition of $D \gencomp_W E$.  The intersection of $U$ 
with $L$ is already the ribbon $\theta_1$, which includes all the cells
of $\south{E}$ for the copies of $E$ that make up $U$.  Therefore, all 
that remains of $U$ is $|\nwD|$ copies of $\nbody{O}$ and $|\nwD|+1$
copies of $\nbody{W}$.  By Hypotheses~\ref{ite:separation}
and \ref{ite:crucial}, these copies of $\nbody{W}$ and $\nbody{O}$ are
pairwise non-adjacent.\no
Construct the southeast decomposition for each copy of $\nbody{W}$ and
$\nbody{O}$.
Label the new ribbons
$\theta_{r+1}, \ldots, \theta_N$ from southwest to northeast, according to 
their southwest endpoints.  Then $\theta_1, \theta_{r+1}, \theta_{r+2}, \ldots, \theta_{N}$
is an outside ribbon decomposition for $U$.  
Therefore, we finally have an outside
ribbon decomposition $\theta_1, \ldots, \theta_N$ for $D \gencomp_W E$.
Let $\Theta$ denote the cutting strip of this decomposition.

\subsection{Applying Sylvester's Determinantal Identity in Case~\ref{ite:staircase}}
With the ribbon decomposition defined in the previous subsection, we now move to the second part of 
the proof of the case when $E=\rightup{W}{O}$.  Exactly 
as in Case~\ref{ite:amalgcomp}\ok we apply Theorem~\ref{thm:sylvester},\ok letting
$M$ be the Hamel-Goulden matrix for the ribbon decomposition $\theta_1, \ldots, \theta_N$ 
and
$S=\{r+1, \ldots, N\}$.  For exactly the same reason as before, the
left-hand side of \eqref{equ:sylvester} becomes
\begin{equation}\label{equ:lhs2}
s_{D \gencomp_W E} \left(\left(s_{\nbody{W}}\right)^{|\nwD|+1}\left(s_{\nbody{O}}\right)^{|\nwD|}\right)^{r-1}.
\end{equation}

We next wish to evaluate $M[S\cup\{i\}, S\cup\{j\}]$ for $i,j \in \{1,2,\ldots,r\}$.  
The top-left entry
of this submatrix is $s_{\theta_i \# \theta_j}$.  Again, we must first show that
if $\theta_i \# \theta_j$ is undefined, then $\det M[S\cup\{i\}, S\cup\{j\}] = 0$.
We omit the proof of this since it is identical to the proof 
for Case~\ref{ite:amalgcomp}.

We next wish to show that a skew diagram 
$F$ with the following property exists:
when using the cutting strip $\Theta$, 
the ribbons of the resulting decomposition of $F$ correspond to the same portions of  
$\Theta$ as $\theta_i \# \theta_j, \theta_{r+1}, \theta_{r+2}, \ldots, \theta_N$. 
Then\ok \linebreak $M[S\cup\{i\}, S\cup\{j\}]$ for $i,j \in \{1,2,\ldots,r\}$ will be the Hamel-Goulden matrix for $F$.

Since $\theta_i \# \theta_j$ is contained in $\theta_1$, and 
$\theta_1, \theta_{r+1}, \theta_{r+2}, \ldots, \theta_N$ are the ribbon decomposition for $U$, 
we wish to obtain $F$ as a subset of $U$.  Our candidate $F$ is much simpler to define than in 
Case~\ref{ite:amalgcomp}: in $U$, just replace $\theta_1$ with the subribbon of $\theta_1$
corresponding to $\theta_i \# \theta_j$.  
In our running example, taking $i=2$ and $j=3$ gives the skew diagram
\[
\begin{matrix}
&&&&&&&&&&&& 12 & 12 \\
&&&&&&&&&&&& \\
&&&&&&&&&&& 11 \\
&&&&&&&&& 10 & 10 \\
&&&&&&&&& 1  \\
&&&&&&&& 9 & 1 \\
&&&&&& 8 & 8 & 1 & 1 \\
&&&&&& 1 & 1 & 1  \\
&&&&& 7 & 1 & \\
&&& 6 & 6 & 1 & 1 \\
&&& 1 & 1 & 1 \\
&& 5 \\
4 & 4
\end{matrix}.
\]

Define $\nbody{E}$, $\nbody{W}\nbody{O}$ and 
$\nbody{O}\nbody{W}$ as in Case~\ref{ite:amalgcomp}.
Writing $\Phi$ for $\phi_1$, the cutting strip of the northwest decomposition of $D$, 
we know that the contribution of $\Phi$ to $D \gencomp_W E$ is $U = E^{\amalg_W |\nwD|}$.
If $\Phi'$ is a portion of $\Phi$ and $E' \subseteq E$, we will let $\Phi' \gencomp_W E'$
denote the subset of $U$ corresponding to $\Phi'$ and $E'$ in the natural way.
Our candidate $F$ can thus also be written as
\[
\left(\Phi[p(\phi_j), q(\phi_i)] \gencomp_W \south{E}\right) \cup
\left(\Phi \gencomp_W \nbody{E}\right) .
\]
This is identical to \eqref{equ:byconstruction}.  Proceeding exactly as in 
Case~\ref{ite:amalgcomp},
we deduce that
\[
\det M[S \cup \{i\}, S \cup \{j\}] = 
s_{(\phi_i \# \phi_j)\gencomp_W E} \left(s_{\nbody{W}}s_{\nbody{O}}\right)^{|\nwD|-q(\phi_i)+p(\phi_j)-1} .
\]
At this point, we diverge from the proof of 
Case~\ref{ite:amalgcomp}.
We still, however,  wish to plug our calculations into \eqref{equ:sylvester}.
In $\det \syl{M}{S}$, we can factor out the powers of $s_{\nbody{W}}s_{\nbody{O}}$ and
divide the result into \eqref{equ:lhs2} as before, to obtain 
\[
s_{D \gencomp_W E} \left(s_{\nbody{W}}\right)^{|\upbody{D}|} \left(s_{\nbody{O}}\right)^{|\nwbody{D}|}
= \det H,\ok
\]
where
\[
H = \left(s_{(\phi_i \# \phi_j)\gencomp_W E}\right)_{i,j=1}^r .
\]
The skew diagram of the $(i,j)$-entry of $H$ is 
$(\phi_i \# \phi_j)\gencomp_W E$, which equals
\[
E^{\amalg_W \left(q(\phi_i)-p(\phi_j)+1\right)}.
\]
Recall that
the ribbons $\phi_1, \ldots, \phi_r$ are those of the enhanced northwest decomposition
of $D$, labelled so that $q(\phi_1) \geq q(\phi_2) \geq \cdots \geq q(\phi_r)$.
We wish to relate these \new{ribbons} to the rows of $D$.  
In fact, since we know that $\det H \neq 0$, we have
$q(\phi_1) > q(\phi_2) > \cdots > q(\phi_r)$.
However, we
certainly need not have $p(\phi_1) > p(\phi_2) > \cdots > p(\phi_r)$.  
To solve this problem, apply a permutation $\sigma^{-1}$ to the columns of
$H$ to obtain a matrix $H'$.  The permutation $\sigma$ is chosen so that
$p(\phi_{\sigma(1)}) > p(\phi_{\sigma(2)}) > \cdots > p(\phi_{\sigma(r)})$.
Then the skew diagram of the
$(i,j)$-entry of $H'$
is 
\[
E^{\amalg_W \left(q(\phi_i)-p(\phi_{\sigma(j)})+1\right)}.
\]
We will have $\det H = (-1)^{\inv{\sigma}} \det H'$, 
where $\inv{\sigma}$ denotes the number of inversions of $\sigma$, 
explaining the appearance of the $\pm$ sign
in Theorem~\ref{thm:gencomp}.  
This sign is completely determined by the enhanced northwest decomposition
of $D$, and so will \new{be} the same for any $E$ of the form $\rightup{W}{O}$.  

To interpret $q(\phi_i)-p(\phi_j)+1$, we need to 
examine the enhanced northwest decomposition of $D$.
The reader may wish to refer to \eqref{equ:enhanced}.
Clearly, the content of the northeasternmost cell of $D$ is $q(\phi_1)$.  Now suppose
$d$ is one position south of $d'$ on $\se{D}$.  
If $d$ and $d'$ are on different ribbons of the northwest decomposition of $D$,
then $d$ is clearly the northeast endpoint of $\phi_i$ for some $i$.
Otherwise $d$ and $d'$ are on the 
same ribbon of the northwest decomposition, in which case the content of $d$ is
$q(\phi_i)$ for some empty ribbon $\phi_i$.  Therefore, 
the content of the easternmost cell on every row of $D$ equals $q(\phi_i)$ for some $i$.
Hence, \new{for $i=1,\ldots,r$, $q(\phi_i)$} must be
the content of the easternmost cell in the $i$th row of $D$.  If we write $D=\lambda/\mu$,
with $\lambda$ chosen so that $\lambda_1$ and the length of $\lambda$ are minimal, then
$q(\phi_i) = \lambda_i -i$.  

Similarly, the content of the southwesternmost
cell of $D$ is clearly $p(\phi_{\sigma(r)})$.  If $d$ is one position south of $d'$ on 
$\nwD$, then $d$ and $d'$ are obviously on the same ribbon, and hence 
every cell on $d$'s diagonal goes north, as defined in Subsection \ref{subsec:HGdets}. 
This means that the cells $e$ and $e'$ on $\se{D}$
with the same contents as $d$ and $d'$ respectively satisfy either:
\begin{itemize}
\item $e$ appears one position south of $e'$ on the same ribbon of the northwest decomposition; or
\item $e'$ is the southwest endpoint of a ribbon of the northwest decomposition
of $D$.
\end{itemize}
In either of the two cases, the content of $e'$, and hence the content of $d'$, equals
{$p(\phi_{\sigma(j)})$} for some $j=1, \ldots, r-1$.  Since there are $r-1$ such $d'$, we
conclude that, for $j=1,\ldots,r$, $p(\phi_{\sigma(j)})$ is the content of the westernmost 
cell in the $j$th row of $D$, i.e.\
$p(\phi_{\sigma(j)}) = \mu_j+1-j$.  

Therefore, $q(\phi_i)-p(\phi_{\sigma(j)})+1 = \lambda_i - \mu_j - i + j$, which is the $(i,j)$-entry of the 
Jacobi-Trudi decomposition matrix for $D$.  Therefore, by Lemma~\ref{lem:importantmap},
$\det H' = s_D \gencomp_W s_E$, and so
\[
s_{D \gencomp_W E} \left(s_{\nbody{W}}\right)^{|\upbody{D}|} \left(s_{\nbody{O}}\right)^{|\nwbody{D}|}
= \pm \left( s_D \gencomp_W s_E \right),
\]
where the sign on the right-hand side is the same for all $E=\rightup{W}{O}$.

\section{Concluding remarks}\label{sec:conclusion}

%

We wish to conclude by making a remark about Conjecture~\ref{con:gencomp} and by introducing two 
further conjectures.

\subsection{Removing Hypothesis\no \ref{ite:final}}

As noted in Conjecture~\ref{con:gencomp}, 
we do not believe that Hypothesis~\ref{ite:final} \new{(see Hypothesis~\ref{hyp:final})} is necessary for Theorem~\ref{thm:gencomp}
to hold.  To prove the first assertion of the conjecture, we need to consider
skew diagrams $E$ such as
\[
E\ =\ \begin{matrix} & & \ti & w \\ & \ti & \ti & w \\ w & \ti \\ w & \ti \end{matrix}.
\]
Using \cite{BucSoftware}, we can check that Conjecture~\ref{con:gencomp} still holds if
$D = \begin{smallmatrix} \ti & \ti \\ \ti \end{smallmatrix}$ or if we put $D^*$ in place
of $D$.  Since $s_D = s_{D^*}$, we conclude that 
$s_{D \gencomp_W E} = s_{D^* \gencomp_W E}$, i.e.
\begin{equation}\label{equ:unexplained}
\begin{matrix}
&&&&&&& \ti & \ti \\
&&&&&& \ti & \ti & \ti \\
&&&& \ti & \ti & \ti \\
&&& \ti & \ti & \ti & \ti \\
&& \ti & \ti \\
&& \ti & \ti \\
& \ti & \ti & \ti \\
\ti & \ti \\
\ti & \ti
\end{matrix}
\sim
\begin{matrix}
&&&&&&& \ti & \ti \\
&&&&&& \ti & \ti & \ti \\
&&&&& \ti & \ti \\
&&&&& \ti & \ti \\
&&&& \ti & \ti & \ti \\
&& \ti & \ti & \ti \\
& \ti & \ti & \ti & \ti \\
\ti & \ti \\
\ti & \ti
\end{matrix}.
\end{equation}
{Since} $E$ does not satisfy Hypothesis~\ref{ite:final}, this equivalence does not
follow from Theorem~\ref{thm:equivs}.  
On the other hand, 
skew-equivalences such as these are explained by Conjecture~\ref{con:gencomp}.
However, since we have been unable to construct outside
ribbon decompositions of the skew diagrams in \eqref{equ:unexplained} that are
amenable to our current techniques, it seems
that some new ideas will be necessary in order to prove Conjecture~\ref{con:gencomp}.

\subsection{Skew diagrams equivalent to their transpose}

It turns out in practice that there are many skew-equivalences of 
the form $F \sim F^t$.  The following result gives an explanation for this.  

\begin{proposition}\label{pro:transpose}
Suppose $E=WOW$ satisfies Hypotheses~\ref{ite:maximal}  -- \ref{ite:crucial} with $E^t=E$,
$W^t=W$ and $W \neq \emptyset$. Then for any skew diagram $D$, 
\[
\left( D \gencomp_W E \right)^t \sim D \gencomp_W E .
\]
\end{proposition}

\begin{proof}
Since $E$ and $W$ are self-transpose, $E$ must be of the form $\rightup{W}{O}$ or
$\upright{W}{O}$.  In particular, $E=WOW$ satisfies \new{Hypothesis~\ref{ite:final} (Hypothesis~\ref{hyp:final})} trivially.
The result now follows from Lemma~\ref{lem:gencomp}(\ref{ite:transpose}) and
Theorem~\ref{thm:equivs}.
\end{proof}

Certainly, if $F=F^t$ then $F \sim F^t$.  
We conjecture that the appropriate converse to Proposition~\ref{pro:transpose} is 
also true.

\begin{conjecture}
Suppose a skew diagram $F$ has the property that $F \sim F^t$, with $F \neq F^t$.  
Then there exists
a skew diagram $E=WOW$ satisfying 
Hypotheses~\ref{ite:maximal}  -- \ref{ite:crucial} and a skew diagram $D$ such that
$F = D \gencomp_W E$ with $E^t=E$, $W^t=W$ and $W \neq \emptyset$.
\end{conjecture}

\subsection{Necessary and sufficient conditions for skew-equivalence}

The strong-est\ok result of \cite{BTV} gives necessary and sufficient conditions for
two ribbons to be skew-equivalent.  
The overarching goal of \cite{RSV07} and the current paper has been to make
progress towards extending this result to general skew diagrams.  We are now in a position
to conjecture such necessary and sufficient conditions.  

First, let us state 
a result that follows by induction from Theorem~\ref{thm:equivs}.

\begin{theorem}\label{thm:compright}
Suppose we have skew diagrams $E_1, E_2, \ldots, E_r$ such that for $i=2,\ldots, r$,
$E_i = W_i O_i W_i$ satisfies Hypotheses~\ref{ite:maximal} -- \ref{ite:final}.  
Let $E'_1$ denote either $E_1$ or $E_1^*$, and for each $i=2,\ldots,r$, let 
$E'_i$ and $W'_i$ denote
either $E_i$ and $W_i$, or $E_i^*$ and $W_i^*$.  Then
\[
\nnew{\left(\cdots\left(\left(}E_1 \gencomp_{W_2} E_2\right) \gencomp_{W_3} E_3 \right) \cdots 
\right) \gencomp_{W_r} E_r \ \sim\
\nnew{\left(\cdots\left(\left(}E'_1 \gencomp_{W'_2} E'_2\right) \gencomp_{W'_3} E'_3 \right) \cdots 
\right) \gencomp_{W'_r} E'_r.
\]
\end{theorem}

Next, let us recall 
Theorem 4.1 from \cite{BTV}
in our notation, where it was also shown that the $\gencomp_\emptyset$ operation is associative
when applied to ribbons.

\begin{theorem}\label{thm:hdl}
\cite[Theorem~4.1]{BTV}\pub Two ribbons $\alpha$ and $\beta$ satisfy $\alpha \sim \beta$ if and only if, 
for some $r$, 
\[
\alpha = \alpha_1 \gencomp_\emptyset \alpha_2 \gencomp_\emptyset \cdots \gencomp_\emptyset
\alpha_r
\mbox{\ \ and\ \ }
\beta = \beta_1 \gencomp_\emptyset \beta_2 \gencomp_\emptyset \cdots \gencomp_\emptyset
\beta_r, 
\]
where, for each $i$, $\alpha_i$ and $\beta_i$ are ribbons with either
$\beta_i = \alpha_i$ or $\beta_i = \alpha_i^*$. The skew-equivalence class of $\alpha$ will contain $2^r$ elements, where $r$ is the number of factors $\alpha_i$ in the irreducible factorisation of $\alpha$ such that $\alpha _i\neq \alpha _i^*$.
\end{theorem}

It transpires that the concept of irreducible factorisation of \cite{BTV} can be extended to arbitrary skew diagrams.

\begin{definition}
Given a \emph{factorisation} of a skew diagram $F=D\gencomp _W E$, where $E=WOW$ satisfies Hypotheses~\ref{ite:maximal} -- \ref{ite:crucial}, we say that the factorisation is \emph{trivial} if 
the factorisation is any one of the following:
\begin{enumerate}
\renewcommand{\theenumi}{\roman{enumi}}
\item $(1) \gencomp _W F$;
\item $F \gencomp _\emptyset (1)$;
\item $\emptyset \gencomp _F E$.
\end{enumerate}
We say the factorisation is \emph{minimal} if it is non-trivial and, among all factorisations
of $F$,  $W$ and then $E$ occupies the minimum number of diagonals.

A factorisation $\nnew{\left(\cdots\left(\left(}E_1 \gencomp_{W_2} E_2\right) \gencomp_{W_3} E_3 \right) \cdots 
\right) \gencomp_{W_r} E_r $  is called \emph{irreducible} if:
\begin{itemize}
\item $E_1$ only admits trivial factorisations;
\item for $i=2,\ldots ,r$ we have $E_i = W_i O_i W_i$ satisfies Hypotheses~\ref{ite:maximal} -- \ref{ite:crucial};
\item each factorisation $D_i \gencomp_{W_i} E_i$ is minimal, where
$$D_i = \nnew{\left(\cdots\left(\left(}E_1 \gencomp_{W_2} E_2\right) \gencomp_{W_3} E_3 \right) \cdots 
\right) \gencomp_{W_{i-1}} E_{i-1}.$$
\end{itemize}
\end{definition}

\begin{remark}\label{rem:notunique}
If $F$ is a ribbon,\ok then we can prove that the irreducible factorisation of $F$ is unique by a proof similar to \cite[Theorem~3.6]{BTV}. Unfortunately, uniqueness does not hold in general, for\ok if $F$ is the left-hand skew diagram of \eqref{equ:unexplained},\ok then it can be irreducibly factored into
$$\begin{matrix}\ti&\ti\\
\ti\end{matrix}\quad \gencomp _{\begin{smallmatrix}\ti\\\ti\end{smallmatrix}}\quad \begin{matrix} & & \ti & \ti \\ & \ti & \ti & \ti \\ \ti & \ti \\ \ti & \ti \end{matrix}.$$However, since $F=F^t$, it can also be factored into
$$\begin{matrix}&\ti\\
\ti&\ti\end{matrix}\quad \gencomp _{\begin{smallmatrix}\ti &\ti\end{smallmatrix}}\quad \begin{matrix} & & \ti & \ti \\ & \ti & \ti & \ti \\ \ti & \ti \\ \ti & \ti \end{matrix}$$by Lemma~\ref{lem:gencomp}(\ref{ite:transpose}). Nonetheless, observe that the number of factors $E_i$ for which $E _i\neq E _i^*$ is the same for both factorisations.
\end{remark}

We now state our main conjecture, of which Theorem~\ref{thm:hdl} implies a very special case.

\begin{conjecture}\label{con:iff}
Two skew diagrams $E$ and $E'$ satisfy $E \sim E'$ if and only if, for some $r$, 
\begin{eqnarray*}
E & = & \nnew{\left(\cdots\left(\left(}E_1 \gencomp_{W_2} E_2\right) \gencomp_{W_3} E_3 \right) \cdots 
\right) \gencomp_{W_r} E_r \mbox{\ \ and} \\
E' & = & \nnew{\left(\cdots\left(\left(}E'_1 \gencomp_{W'_2} E'_2\right) \gencomp_{W'_3} E'_3 \right) \cdots 
\right) \gencomp_{W'_r} E'_r,
\end{eqnarray*}
where 
\begin{itemize}
\item $E_1, E_2, \ldots, E_r$ are skew diagrams;
\item for $i=2,\ldots, r$,
$E_i = W_i O_i W_i$ satisfies Hypotheses~\ref{ite:maximal} -- \ref{ite:crucial};
\item $E'_1=E_1$ or $E'_1=E_1^*$;
\item for $i=2,\ldots,r$, $E'_i$ and $W'_i$ denote
either $E_i$ and $W_i$, or $E_i^*$ and $W_i^*$.
\end{itemize} The skew-equivalence class of $E$ will contain $2^r$ elements, where $r$ is the number of factors $E_i$ in any irreducible factorisation of $E$ such that $E _i\neq E _i^*$.
\end{conjecture}

\begin{remark} First,\ok note that it was conjectured \cite[Conjecture~9.1]{RSV07} that skew-equivalence classes have size a power of 2, and with our construction we can now predict precisely which power. Secondly, observe that\new{ }our definition of irreducible factorisation for ribbons differs from that in \cite{BTV}.  As an example, the  diagram $(4)$ is irreducible under the definition in \cite{BTV} but is reducible to $(2)\gencomp _\emptyset (2)$ under our definition. However, 
the powers of 2 in Theorem~\ref{thm:hdl} and Conjecture~\ref{con:iff} are indeed the same. The condition of minimising $W$ first in the irreducible factorisation for skew diagrams ensures that we can naturally pair up terms in the expressions arising from both statements,\ok and those terms that do not pair up will not contribute to the power of 2.

For an example
of a  skew-equivalence class of size greater than 4, we see that in Remark~\ref{rem:staircase},
$D \gencomp_W E$ is actually equal to $(E \gencomp_W E) \gencomp_W E$, where
$E=(2,1)$ and $W=(1)$.  One can check that the 8 skew diagrams of the form
$(E' \gencomp_W E'') \gencomp_W E'''$, with $E'$, $E''$ and $E'''$ each equal to 
$E$ or $E^*$, are all different.  
\end{remark}

To conclude we present evidence in favour of Conjecture~\ref{con:iff}. Observe that the 
``if'' direction of Conjecture~\ref{con:iff}
would follow from Conjecture~\ref{con:gencomp} in the
same way that Theorem~\ref{thm:compright} follows from Theorem~\ref{thm:gencomp}. 
The only difference is that Hypothesis~\ref{ite:final} is absent in the conjectures.
To support both the converse direction and the skew-equivalence class sizes, we have verified that the conjecture holds for all  skew diagrams with at most 20 cells. 

\bibliography{skeweq}
\bibliographystyle{plain}
%
%

\end{document}